

        \documentstyle{amsppt}
        \NoBlackBoxes
        \magnification=\magstep1
\define\longerrightarrow{\DOTSB\relbar\joinrel\relbar\joinrel
\relbar\joinrel\rightarrow}

       \topmatter
        \title  Infinite Grassmannians and Moduli Spaces of
       $G$--bundles
       \endtitle
        \author  Shrawan Kumar
      \endauthor
        \endtopmatter
        \document

    \vskip4ex

\hyphenation{asso-cia-ted gene-ra-li-sed}

 \subheading{Introduction}

These are notes for my eight lectures given at the C.I.M.E. session on ``Vector
bundles on curves. New directions"  held at Cetraro (Italy) in June 1995. The
work presented here was done in collaboration with M.S. Narasimhan and A.
Ramanathan and appeared in \cite{KNR}. These notes differ from \cite{KNR}
in that we have added three appendices (A)-(C) containing basic definitions and
results (we need) on ind-varieties, Kac-Moody Lie algebras, the associated
groups and their flag varieties. We also have modified the proof (given in \S7)
of basic extension result (Proposition 6.5), and we hope that it is more
transparent than the one given in \cite{KNR, \S7}. We now describe the main
result of this note.

  Let $C$ be a smooth projective irreducible algebraic
 curve over $\Bbb C$ of any genus and $G$ a connected  simply-connected
  simple affine algebraic group over $\Bbb C$.  In this note we
 elucidate the relationship between
 \roster
 \item the space of vacua (``conformal blocks'') defined
 in Conformal Field Theory, using an integrable highest weight  representation
of the
 affine Kac-Moody  algebra associated to $G$ and
 \item the space of regular sections (``generalised theta functions'')
 of a line bundle on the moduli space $\frak M$ of semistable
 principal $G$-bundles on $C$.
 \endroster

 \vskip1ex

 Fix a point $p$ in $C$ and let $\hat{\Cal O}_p$ (resp.
 $\hat{\text
{{\sl k}}}_p$) be the completion of the local ring
 $\Cal O_p$ of $C$ at $p$ (resp. the quotient field
  of $\hat{\Cal O}_p$).  Let $\Cal G := G(\hat{\text
{{\sl k}}}_p$)  (the $\hat{\text
{{\sl k}}}_p$-rational points of the algebraic group $G$)
 be the loop group of $G$ and let $\Cal P := G(\hat{\Cal O}_p)$
 be the standard maximal parahoric subgroup of $\Cal G$.  Then the generalised
flag variety $X :=\Cal
 G/\Cal P$ is an inductive limit of projective varieties,
 in fact of generalised Schubert varieties.  One has a
 basic homogeneous  line bundle $\frak L(\chi_o)$ on $X$ (cf. \S C.6), and
 the Picard group Pic$(X)
$ is isomorphic to
 $\Bbb Z$ which is generated by $\frak L(\chi_o)$
 (Proposition C.13).  There is a
 central extension $\tilde{\Cal G}$ of $\Cal G$ by the multiplicative group
$\Bbb C^{\ast}$ (cf. \S C.4),
which acts on the line bundle $\frak L(\chi_o)$.  By an analogue of the
Borel-Weil theorem proved in the Kac-Moody setting by Kumar (and also by
Mathieu), the space $H^0(X,\frak L(d\chi_o)) $ of the regular sections of the
line bundle
 $\frak L(d\chi_o)$ :=$\frak L(\chi_o)^{\otimes d} $ (for any $d \geq 0$)
 is canonically isomorphic with the full vector space dual $L(\Bbb C,d)^\ast$
of the integrable
 highest weight (irreducible) module $L(\Bbb C,d)$ (with central charge $d$)
of the affine Kac-Moody Lie algebra $\tilde{\frak g}$
 (cf. \S A.2).
 \vskip1ex

 Using the fact that any principal $G$-bundle on
 $C\setminus p$ is trivial (Proposition 1.3), one sees
 easily that the set of isomorphism classes of principal
 $G$-bundles on $C$ is in bijective correspondence with
 the double coset  space $\Gamma \backslash \Cal G/\Cal
 P$, where $\Gamma :=$ Mor$(C\setminus p,G)$ is the subgroup
of $\Cal G$ consisting of all the algebraic morphisms  $C\setminus p
\to G$.  Moreover,
 $X$ parametrizes an algebraic family
$\Cal U$
 of principal
 $G$-bundles on $C$ (cf. Proposition 2.8). As an interesting  byproduct of this
parametrization and rationality of the generalised Schubert varieties $X_{\frak
w}$, we obtain that the moduli space $\frak M $ of semistable principal
$G$-bundles on $C$
is  a unirational variety (cf. Corollary 6.3). Now, given  a finite dimensional
 representation $V$ of $G$, let ${\Cal U}(V$) be the family of
 associated vector bundles on $C$ parametrized by $X$.
 We have then the determinant line bundle $\text{Det}(\Cal U(V$)) on $X$,
 defined as the dual of the determinant of the cohomology of the
 family $\Cal U(V)$
  of vector bundles on $C$ (cf. \S 3.7). As we mentioned above,
Pic ($X$) is freely generated by the homogeneous line bundle
 $\frak L $($\chi_o$)  on $X$, in particular, there exists a unique
integer $m_V$ (depending on the choice of the representation $V$)
such that
 $\text{Det}(\Cal U (V)) \simeq  \frak L(m_V\chi_o)$. We
determine this number explicitly in Theorem (5.4), the proof of which
makes use of Riemann-Roch theorem. It is shown  that the number
$m_V$  coincides with the Dynkin index of the representation $V$.
For example,  if we take $V$ to be the adjoint representation of
$G$, then $m_V$ =$ 2 \times$ dual Coxeter number of $G$ (cf. Lemma 5.2 and
Remark 5.3). The number $m_V$ is also expressed in terms of the induced map at
the third homotopy group level $\pi_3(G) \to \pi_3(SL(V))$ (cf. Corollary 5.6).
 \vskip1ex

 The action of  $\Gamma $  on $X$ via left multiplication lifts to an action on
the line bundle $\frak L(m_V\chi_o)$ (cf. \S 2.7). Suggested by Conformal Field
Theory, we consider the
space $H^0(X,\frak L(dm_V\chi_o))^
 {\Gamma}$ of $\Gamma$-invariant regular sections of the   line bundle
 $\frak L(dm_V\chi_o)$ (for any $d \geq 0$).  This space of invariants is
 called the space of vacua. More precisely, in Conformal Field Theory,  the
space of vacua is defined to be the space of invariants of the Lie algebra
$\frak g \otimes R$ in $ L(\Bbb C,d)^\ast$, where $R$ is the ring of regular
functions on the affine curve $C\setminus p$ and  $\frak g$ is the Lie algebra
of the group $G$.  We have (by Proposition 6.7) $ [L(\Bbb C,dm_V)^\ast]^\Gamma
=  [L(\Bbb C,dm_V)^\ast]^{\frak g \otimes R}$ and, as
already mentioned above,  $H^0(X,\frak L(dm_V\chi_o)) \simeq L(\Bbb
C,dm_V)^\ast.$  The main result of this note (Theorem 6.6) asserts  that  (for
any $d \geq 0$) the space $ H^0(\frak M, \Theta (V)^{\otimes d})$ of the
regular sections
of the $d$-th power of the $\Theta$-bundle $\Theta(V)$ (cf. \S3.7) on the
moduli space $\frak M$ is  isomorphic with the space
of vacua  $ [L(\Bbb C,dm_V)^\ast]^\Gamma =[L(\Bbb C,dm_V)^\ast]^{\frak g
\otimes R} $.  This is  the connection, alluded to in the
beginning of the introduction, between the space of vacua and the space
of generalised theta functions.
   This
 result has also independently been obtained   by Faltings \cite{Fa} and in the
case of $G = SL_N $ by Beauville-Laszlo \cite{BL}, both  by different methods .
\vskip1ex

 We make crucial use of a `descent' lemma (cf. Proposition 4.1), and an
extension result (cf. Proposition 6.5)
in the proof of Theorem (6.6). The proof of Proposition (6.5) is given in \S7,
and relies on the explicit GIT construction of the moduli space of vector
bundles.

Our Theorem (6.6) can be generalised to the situation where the
curve $C$ has $n$ marked points $\{p_1, \dots, p_n\}$ together with finite
dimensional  $G$-modules $\{V_1, \dots, V_n \}$ attached to them
respectively, by bringing in moduli space of parabolic
$G$-bundles on $C$.

   A purely algebro-geometric
study (which does not use loop groups) of generalised theta
fuctions on the moduli space of (parabolic) rank two torsion-free
sheaves on a nodal curve is made by Narasimhan-Ramadas \cite{NRa}.
A factorization theorem  and a vanishing theorem for the
theta line bundle are proved there. In addition, several other mathematicians
(A. Bertram, S. Bradlow, S. Chang, G. Daskalopoulos, B. Geemen, E. Previato, A.
Szenes, M. Thaddeus, R. Wentworth, D. Zagier, $\cdots$ ) and physicists have
studied
the space of generalised theta functions (from different view points) in the
case when $G = SL(2)$, in the last few years.

\vskip1ex

The organization of the note is as follows:

Apart from introducing some  notation in \S 1, we  realize the affine flag
variety $X$ as a parameter set for $G$-bundles. In section (2) we  prove that
$X$ supports  an algebraic family of G-bundles on the curve
$C$ (cf. Proposition 2.8). We also realize the group $\Gamma$ as an ind-group,
calculate its Lie algebra, and prove its splitting in this section. Section (3)
 is devoted to recalling some basic definitions and results on the moduli space
of semistable $G$-bundles, including the definition of the determinant line
bundle and
the $\Theta$-bundle on the moduli space. We prove a
curious result (cf. Proposition 4.1) on algebraic descent in \S 4. Section (5)
is devoted to identifying the determinant line bundle on $X$
with a suitable power of the basic homogeneous line bundle on $X$.  Section (6)
contains the statement and the proof of the main result (Theorem 6.6).
Finally in Section (7) we prove the basic extension result (Proposition 6.5),
using Geometric Invariant Theory.
Appendix (A) is devoted to recalling the definition of affine Kac-Moody Lie
algebras and its representations. Appendix (B) is an introduction to
ind-varieties and ind-groups. Finally in appendix (C), we recall the basic
theory of affine Kac-Moody groups and their flag varieties.
\vskip9ex

   \subheading{ 1. Affine flag variety as parameter set for
       $G$--bundles}
  \vskip2ex
 \flushpar ({\bf 1.1}) {\it  Notation.}  Throughout the note
       $k$ denotes an algebraically closed field of char. 0.
     By a scheme we will mean a scheme over $k$. Let us fix a
      smooth irreducible  projective curve $C$ over
       $k$, and a  point $p\in C$.  Let $C^{\ast}$ denote
    the open set $C\setminus p$. We also
       fix an affine algebraic connected simply-connected simple  group $G$
       over $k$.

        For any $k$-algebra $A$, by $G(A)$ we mean the
       $A$--rational points of the algebraic group $G$.  We fix
       the following notation to be used throughout the note:
        \roster
        \item"${}$" $\Cal G = \Cal G_{\frak T} =
        G(\hat{\text{\text{{\sl k}}}}_p)$,
        \item"${}$" $\Cal P = \Cal P_{\frak T} = G(\hat{\Cal
       O}_p)$, and
        \item"${}$" $\Gamma = \Gamma_{\frak T} =
       G(k[C^{\ast}])$,
        \endroster
        where $\hat{\Cal O}_p$ is the completion of
       the local ring $\Cal O_p$ of $C$ at $p$,
       $\hat{\text{{\sl k}}}_p$ is the quotient field of
       $\hat{\Cal O}_p$, $k[C^{\ast}]$ is the ring of
       regular functions on the affine curve $C^{\ast}$
       (which can canonically be viewed as a subring of
       $\hat{\text{{\sl k}}}_p$), and ${\frak T}$ is the triple
       $(G,C,p)$.

        We recall the following \vskip.5ex
        \flushpar ({\bf 1.2}) {\it   Definition}.  Let $H$ be any (not
necessarily reductive) affine
   algebraic group over $k$ . By a {\it principal
       $H$-bundle} (for short $H${\it -bundle}) on an
       algebraic variety $X$, we mean an algebraic variety $E$
       on which $H$ acts algebraically from the right and a
       $H$-equivariant morphism $\pi :E  \to X$ (where $H$ acts
       trivially on $X$), such that $\pi $ is isotrivial (i.e.
       locally trivial in the \'etale topology).
        \vskip.5ex

        Let $H$ act algebraically on a quasi-projective variety
       $F$ from the left.  We can then form the {\it associated
       bundle with fiber} $F$, denoted by $E(F)$.  Recall that
       $E(F)$ is the quotient of $E\times F$ under the
       $H$-action given by $g(e,f)= (eg^{-1},gf)$, for $g\in
       H$, $e\in E$ and $f\in F$.

        {\it Reduction of structure group of $E$ to a closed
       algebraic subgroup} $K\subset H$ is, by definition, a $K$-bundle
       $E_K$ such that $E_K(H)\approx E$, where $K$ acts on $H$
       by left multiplication.  Reduction of structure group to
       $K$ can canonically be thought of as a section of the
       associated bundle $E(H/K)\to X$.

        Let $\Cal X= \Cal X(H,C)$ denote the set of isomorphism
       classes of $H$-bundles on the base $C$, and $\Cal X_0=
       \Cal X_0({\frak T})\subset \Cal X$ denote the subset
       consisting of those $H$-bundles on $C$ which are
       algebraically trivial restricted to $C^{\ast}$. We recall the following
proposition essentially due to Harder \cite{H$_1$, Satz 3.3 and the remark
following it}.

        \proclaim{ (1.3)  Proposition}  {\it Let $H$ be a
       connected reductive algebraic group over $k$.  Then the
       structure group of a $H$-bundle on a smooth affine
       curve $Y$ can be reduced to the connected component
       $Z^0(H)$ of the centre $Z(H)$ of $H$.

        In particular, if $H$ as above is semi-simple, then any
       $H$-bundle on $Y$ is trivial.}\endproclaim

        \vskip1ex

        The following map is of basic importance for us in this
       note.  This provides a bridge between the moduli space
       of $G$-bundles and the affine (Kac-Moody) flag variety, where $G$ is as
in \S1.1.
        \vskip.5ex
        \flushpar  ({\bf 1.4})  {\it  Definition} (of the map $\varphi
       :\Cal G\to \Cal X_0$).   Consider the canonical
       morphisms $i_1:\,\text{Spec}\,(\hat{\Cal O}_p)\to C$ and
       $i_2: C^{\ast} \hookrightarrow C$.  The morphisms
       $i_1$ and $i_2$ together provide a flat cover of $C$.
       Let us take the trivial $G$-bundle on both the schemes
       Spec$\,(\hat{\Cal O}_p)$ and $C^{\ast}$.  The fiber
       product
        $$
        F := \,\text{Spec}\,(\hat{\Cal O}_p) \underset C\to{\times}
        C^{\ast}
        $$
        of $i_1$ and $i_2$ can canonically be identified with
       Spec$\,(\hat{\text{{\sl k}}}_p)$. This identification
       $F \simeq$ Spec $(\hat{\text{{\sl k}}}_p)$ is induced from the
       natural morphisms
       $$
       \matrix
       &&\text{Spec }(\hat{\text{{\sl k}}}_p )\\
       &\swarrow  &&\searrow\\
       &\text{Spec } ( \hat{\Cal O}_p )
       &@VVSV \negthickspace \negthickspace \negthickspace
       \negthickspace \negthickspace \negthickspace C^\ast\\
       &\nwarrow  &&\nearrow\\
       &&F
       \endmatrix
       $$

       \par By a \lq\lq glueing\rq\rq  lemma of Grothendieck \cite{Mi, Part I,
Theorem 2.23, p.19}, to give a
       $G$-bundle on $C$, it suffices to give an automorphism
       of the trivial $G$-bundle on Spec$\,(\hat{\text{{\sl
       k}}}_p)$, i.e., to give an element of $\Cal G:=
       G(\hat{\text{{\sl k}}}_p)$.   (Observe
       that since we have a flat cover of $C$ by only two
       schemes, the cocycle condition is vacuously satisfied.)
       This is, by definition, the
       map $\varphi : \Cal G\rightarrow \Cal X_0$.

 \proclaim{ (1.5)  Proposition} {\it The map $\varphi$
       (defined above) factors through the double coset space
       to give a bijective map (denoted by)}
       $$
       \bar{\varphi} : \Gamma \backslash \Cal G/\Cal P
       \rightarrow
       \Cal X_0\,.
       $$
       \endproclaim
       (Observe that, by Proposition (1.3), $\Cal X_0=\Cal X$
       since $G$ is assumed to be connected and semi-simple.)

       \demo{Proof}  From the above construction, it is clear
       that for $g,g'\in\Cal G,\ \varphi (g)$ is isomorphic
       with $\varphi (g')$ (written $\varphi (g)\approx \varphi
       (g')$) if and only if there exist two $G$-bundle
       isomorphisms :
       $$
       \gather
       \text{Spec } (\widehat{\Cal O}_p )\times G
       \overset\theta_1\to{\underset\sim\to\longrightarrow}
       \text{Spec } (\widehat{\Cal O}_p )\times G \\
       \searrow \qquad \qquad \qquad\swarrow\\
       \text{Spec } (\widehat{\Cal O}_p)\\
       \intertext{and}
       C^{\ast}\times G
       \overset\theta_2\to{\underset\sim\to\longrightarrow}
       C^{\ast}\times G \\
       \searrow \qquad \qquad\swarrow\\
       C^\ast
       \endgather
       $$
       such that the following diagram is commutative:
       $$\CD
       \text{Spec}\,(\hat{\text{{\sl k}}}_p )\times G
       @>\sideset^{\theta_1}\and\to{\vert}_{\text{Spec}\,
       (\hat{\text{{\sl k}}}_p )}>>
       \text{Spec}\,(\hat{\text{{\sl k}}}_p )\times G\\
       @VVg'V   @VVgV\\
       \text{Spec}\,(\hat{\text{{\sl k}}}_p )\times G
       @>\sideset^{\theta_2}\and\to{\vert}_{\text{Spec}\,
       (\hat{\text{{\sl k}}}_p )}>>
       \text{Spec}\,(\hat{\text{{\sl k}}}_p )\times G\ .
       \endCD \tag"$(\ast)$"
$$
       \par Any $G$-bundle isomorphism $\theta_1$ (resp.
       $\theta_2$) as above is given by an element $h\in\Cal
       P$ (resp. $\gamma \in\Gamma $).  In particular, from
       the commutativity of the above diagram ($\ast$),
       $\varphi (g)\approx \varphi (g')$ if and only if there
       exists $h\in\Cal P$ and $\gamma \in\Gamma $ such that
       $gh=\gamma g'$, i.e., $\gamma ^{-1}gh=g'$.  This shows that
       the map $\varphi$ factors through $\Gamma \backslash
       \Cal G/\Cal P$ to give an injective map $\bar{\varphi}$.
       The surjectivity of $\bar{\varphi}$ follows immediately
       from the definition of $\Cal X_0$, and the fact that any
       $G$-bundle on Spec$\,(\widehat{\Cal O}_p )$ is trivial.
       \enddemo
       \vskip1ex
        \flushpar {\bf (1.6)} {\it   Remark}.
       $\Cal G/\Cal P$ should be thought of
       as a parameter space for $G$-bundles $E$ together
       with a trivialization of
       $E_{\vert_{C^\ast}}$ (cf. Proposition 2.8).

\vskip4ex

\subheading {2. Affine flag variety parametrizing an algebraic family and
realizing $\Gamma $  as an ind-group}
\vskip3ex

Recall the definition of the group $\Gamma \subset \Cal G$ from \S1.1.

\vskip 2mm

 \flushpar
 {\bf (2.1) Lemma.}  {\it The group $\Gamma$ is an ind-group.}
 \vskip 3mm
 \flushpar
 {\bf Proof.}\footnote{I thank R. Hammack for some simplification in my
original argument.} Embed $G \subset SL_N \hookrightarrow M_N$, where
  $M_N = M_N(k)$ is the space of $N \times N$ matrices over $k$.  Take a
 $k-$basis $\{f_1,f_2,f_3, \cdots \}$ of $k [C^\ast]$ (the ring of
 regular functions) such that $\text{ord}_pf_n \leq \text{ord}_p
 f_{n+1}$ for any $n\geq 1$, where $\text{ord}_p f_n$ denotes the order of the
pole of $f_n$ at $p$.
 There is an injective map   $i:\Gamma \hookrightarrow \text{Mor} (C^*,
 M_N)$, where Mor denotes the set of all the morphisms.  The set Mor
 $(C^*,M_N)$ has a filtration Mor$_0 \subseteq \dots \subseteq $ Mor$_n
\subseteq \dots ,  $ where Mor$_n $ is the (finite dimensional)  vector space
 of all those morphisms $\theta:C^* \rightarrow M_N$ such that
 all its matrix entries have poles of order $\leq n$.  Set
 $\Gamma_n=i^{-1} (\text{Mor}_n)$.  Any $\theta= (\theta_{i,j}) \in
\text{Mor}_n
 $ can be written as $\theta_{i,j}=
 {\displaystyle{\sum_{k=1}^{k(n)}}} z_{i,j}^{k} f_k$ (for some
 $k(n))$.  We take  $(z_{i,j}^{k})$ as the coordinates on
 Mor$_n$.  It is  easy to see that $\Gamma_n
 \hookrightarrow \text{Mor}_n$ is given by the vanishing of some
 polynomials in  $(z_{i,j}^{k})$ ,  in particular, $\Gamma_n$ is a closed
 subvariety of the affine space Mor$_n$.  (We put the
 reduced  structure on $\Gamma_n$.)  This gives rise to the
 ind-variety structure on $\Gamma$ as a closed ind-subvariety of
 Mor $(C^*, M_N)$.  It is easy to see  (from the definition of the ind-variety
structure on $\Gamma$)  that
 $\Gamma $ in fact is an ind-group.  Moreover,  this ind-variety structure  on
$\Gamma$ does not depend upon the particular choice of
the embedding  $G \hookrightarrow SL_N$. \qed

\vskip1ex

The following lemma determines the Lie algebra of the  ind-group $\Gamma$.
\proclaim{(2.2) Lemma}  The Lie algebra $\text{Lie}\, \Gamma$ is  isomorphic
with $\frak g
\otimes_k k[C^\ast]$,  where $\frak g$ :=  Lie $G$  and the bracket in  $\frak
g
\otimes k[C^\ast]$ is defined as $[X\otimes p, Y\otimes q] = [X, Y]\otimes pq$,
for $ X, Y \in$  $\frak g$ and $ p,q \in  k[C^\ast]$.  The isomorphism
$\text{Lie} \, \Gamma \simeq  \frak g
\otimes k[C^\ast]$ is obtained  by considering
the differential of the evaluation map at each point of $C^\ast$.
\endproclaim

\demo{Proof} Choose an
embedding  $G \hookrightarrow
SL_N  \subset M_N$ as in the proof of Lemma (2.1).   This gives rise to a
closed immersion
$i :  \Gamma  \hookrightarrow $ Mor $(C^\ast, M_N)$.   In particular, it
induces an injective map  $di :
T_e(\Gamma) = \text{Lie} \,\Gamma \hookrightarrow T_I$(Mor) $\simeq $ Mor at
the Zariski tangent space level (where $I$ is the identity matrix and Mor = Mor
$(C^\ast,M_N)).$ We claim that  $di$ is a Lie algebra homomorphism, if we endow
 Mor $\simeq M_N(R)$
with the standard Lie algebra structure, where $R:= k[C^\ast]$ and $M_N(R)$ is
the  space of $N\times N$ matrices over $R$. To prove this, consider the
following commutative diagram (for any fixed $x \in C^\ast$):

 $$\gather
       T_e(\Gamma)
       \ \ \overset{d i}\to \hookrightarrow \ \ \
     M_N(R)\\
       \downarrow \quad \qquad \quad \downarrow\\
       {\frak g}=T_e(G)\   \hookrightarrow \ \
       M_N(k)\  ,
       \endgather$$
where the vertical maps are induced by the evaluation map $e_x : R \to k $
given by $p \mapsto p(x)$. Since the bottom horizontal map  is a Lie algebra
homomorphism, and
so are the vertical maps, we obtain that $di$ itself is a Lie
algebra homomorphism. It is further clear, from the above commutative diagram,
that the image of $di$ is contained in $\frak g
\otimes R$, where $\frak g$  is identified with its image in $M_N(k)$.

Next, we prove that the image of $di$ contains at least the set
$\frak g
\otimes R$:

Fix any  vector $X \in \frak g \subset M_N(k)$ such that $X$ is a nilpotent
matrix and take $p \in R$, and define
a morphism $\Bbb A^1 \to \Gamma$ by $z \mapsto exp\, (z X\otimes p).$
(Since $X$ is nilpotent, the image is indeed contained in $\Gamma$.) It is easy
to see that the image of the induced map (at the tangent space  level at $0$)
is precisely the space $ k(X\otimes p)$. But since the  nilpotent matrices  $X
\in \frak g$    span $\frak g$, the assertion follows. This completes the proof
of the lemma. \qed
\enddemo

 \vskip1ex

We prove the following curious lemma.
\proclaim{(2.3) Lemma}  Let $X$ be a connected variety over $\Bbb C$.  Then any
regular map $ X\to \Bbb C^\ast$, which is null-homotopic in the topological
category, is a constant.\endproclaim

(Observe that if the singular cohomology $H^1(X,\Bbb Z)=0$, then any continuous
map $ X\to \Bbb C^{\ast}$ is null-homotopic.)
\vskip1ex
\demo{Proof}  Assume, if possible, that there exists a null-homotopic
non-constant regular map $\lambda:X\to\Bbb C^\ast$. Since $\lambda$ is
algebraic,  there exists a number $N>0$ such that the number of irreducible
components of $\lambda^{-1}(z)\leq N$, for all $z\in\Bbb C^\ast$.
  Now we consider the $N'$-sheeted covering
$\pi_{N'}:\Bbb C^\ast\to\Bbb C^\ast(z\mapsto z^{N'})$, for any $N'>N$.  Since
$\lambda$ is null-homotopic,
  there exists a (regular) lift
$\tilde{\lambda} :X\to\Bbb C^\ast$ (cf. \cite{Se$_1$, Proposition 20}), making
the following diagram commutative:
$$
\aligned
{}&{\Bbb C}^\ast \\
\sideset^{\tilde{\lambda}\negthickspace\negthickspace}\and\to
\nearrow
\qquad
& \downarrow \pi_{\negthickspace {}_{N'}}V \\
X\quad  \underset\lambda\to\longrightarrow \quad &{\Bbb C}^\ast ~.
\endaligned
$$
Since $\tilde{\lambda}$ is regular and non-constant, by Chevalley's theorem,
Im $\tilde{\lambda}$ (being a constructible set) misses only finitely many
points of
${\Bbb C}^\ast$.  In particular, there exists a $z_o\in \Bbb C^\ast$ (in fact a
Zariski open set of points) such that $\pi^{-1}_{N'}(z_o)\subset
\text{Im}\,\tilde{\lambda}$. But then the number of irreducible components of
$\lambda^{-1}(z_o) = \tilde{\lambda}^{-1}(\pi^{-1}_{N'})(z_o) \geq N'>N$, a
contradiction to the choice of $N$.  This proves the lemma. \qed
\enddemo
\vskip1ex

We will use the following corollary in the proof of assertion (c) contained in
the proof of Theorem (6.6)
\proclaim{(2.4) Corollary}  Take the base field $k={\Bbb C}$. There does not
exist any non-constant regular map
$\lambda:\Gamma\to\Bbb C^\ast$.
\endproclaim

\vskip1ex

\demo{Proof}    The group $\Gamma$ is connected and
simply-connected, in particular, $H^1(\Gamma,\Bbb Z)=0$, where
$H^\ast(\Gamma,\Bbb Z)$ denotes the singular cohomology of the topological
space $\Gamma$.  (In fact, it is very likely that the space $\Gamma$ is
homotopically  equivalent to the corresponding space $\Gamma_{\text{top}}$
consisting of all the continuous maps  $C^\ast \to G$ under the compact-open
topology. This of course will give the connected and simple-connectedness of
$\Gamma$, by using a result of Thom \cite{GK, Theorem 5.10}. A student of mine
R. Hammack is working to write a proof of this homotopy equivalence by using
some ideas similar to \cite{PS, Proof of Proposition 8.11.6}, albeit in the
algebraic category,  together  with a variant of a result of Hurtubise
\cite{Hur, Theorem 1.1}.) This gives that the map
$\lambda$ is null-homotopic.  By using the above lemma (2.3),
$\lambda$ is constant on each connected component of $\Gamma_n$
(for any $n\geq 0)$ and hence
$\lambda$ itself is constant.
\qed
\enddemo

Recall from Proposition (C.12) that $X=X_{\text{rep}}=X_{\text{lat}}$ is a
projective ind-variety.
\vskip1ex
 \flushpar
 {\bf (2.5) Lemma.} {\it The left multiplication of $\Gamma$ on $X$ is a
morphism } $\varphi :\Gamma \times X \to X$.
 \vskip 2mm
 \flushpar
 {\bf Proof.} We will  consider the
 $X_{\text{lat}}$ description of $X$ (cf. \S C.9).   It is clear that
 $\varphi (\Gamma_n \times \hat{X}_m) \subseteq  \hat{X}_{k
 (n,m)},$ for some $k(n,m)$.  Now from the explicit description
 of the variety structures (on $\Gamma$ and $X_{\text{lat}}$), it is easy to
check that
 $\varphi_{n,m}:= \varphi_{\mid \Gamma_n \times \hat{X}_m}$
 is a morphism.

 This proves the lemma. $\qed$
\vskip1ex

Restrict the central extension (1) of \S (C.4) to get a central
 extension
$$
 1 \rightarrow {\Bbb C}^\ast\underset i\to\rightarrow
 \tilde{\Gamma}
 \underset\pi\to\rightarrow \Gamma \rightarrow 1\, ,
  \tag"(1)"$$
where $\tilde{\Gamma}$ is by definition $\pi^{-1}(\Gamma)$.

 \flushpar \subhead (2.6) Splitting of the central
 extension over $\Gamma$ ($SL_N$ case) \endsubhead  The basic reference for
this subsection is \cite{PS, \S7.7}. We first consider
 the case of $G=SL_N$ and follow the same notation as in
  \S (C.7).  In particular, $\Cal G^o := SL_N(\Bbb
 C((t)))$, $\Cal P^o= SL_N(\Bbb C [[t]])$, $X^o= \Cal
 G^o/ \Cal P^o$, $V=\Bbb C^N$, $V((t)) = V\otimes_{\Bbb
 C}\Bbb C((t))$, and $L_o=V\otimes_{\Bbb C}\Bbb C[[t]]$.  Let $GL(W)$ denote
the group of $\Bbb
 C$-linear isomorphisms of a vector space $W$.

  Define the subgroup $\Cal H$ of $\Cal G^o\times
 GL(L_o)$ by
 $$
 \Cal H = \{ (g,E) \in \Cal G^o\times GL(L_o) :
 g^+ -E:L_o\to L_o \text{ has finite rank}\}
 $$
 where $g=\pmatrix g^+ &\ast\\ \ast &\ast\endpmatrix$
 with respect to the decomposition
 $$
 V((t)) = L_o \oplus (V\underset{\Bbb C}\to{\otimes} t^{-1}\Bbb C[t^{-1}]) .
 \tag1$$

  Let $\Cal N\subset \Cal H$ be the normal subgroup
 defined as $\Cal N= \{ (1,E)\in \Cal H:\det E=1\}$.
 (Observe that since $I-E:L_o\to L_o$ has finite rank, i.e., has
 finite dimensional image, the
 determinant of $E$ is well defined.)

 It is not difficult to see that the projection on the
 first factor gives rise to a central extension:
 $$
 1\to \Bbb C^\ast \to \Cal H/\Cal N \to \Cal G^o \to 1~.
 \tag2$$

 We now give an alternative description of the line
 bundle ${\frak L}(\chi_0)$ on $X^o$ (cf. \S C.6):

 Recall the definition of the set $\Cal F$ and the map
 $\beta :X^o\to \Cal F$ from \S (C.7).  For any $W\in
 \Cal F$, define $S_W$ as the set of $\Bbb C$ -linear
 isomorphisms $\theta: L_o\to W$ such that
 $\pi_1\theta-I:L_o\to L_o$ has finite rank, where $\pi_1:V((t))\to L_o$
 is the projection on the $L_o$ factor with respect to
 the decomposition (1).

 Define a vector space $\Cal V_W$ over $\Bbb C$  with
 basis parametrized by the set $S_W$, i.e., an element of
 $\Cal V_W$ is of the form $\sum_{\theta\in
 S_W}z_{\theta}\theta$, where all but finitely many
 $z_{\theta}$ are zero.  Let  $\Cal V '_W
 \subset \Cal V_W$ be  the subspace spanned by $\{ \theta- \det
 (\theta\theta^{'-1})\theta'\}_{\theta,\theta'\in S_W}$
 and let $\Cal L_W= \Cal V_W/\Cal V '_W$.  (Note that
 $\theta\theta^{'-1}-I$ has finite rank as an
 endomorphism of $L_o$ and hence
 $\det(\theta\theta^{'-1})$ is well defined.)   Then
 $\Cal L_W$ is a 1-dimensional vector space.  Now define
 the line bundle $\Cal L\overset
\eta\to{\rightarrow}\Cal
 F$, where $\eta^{-1}(W)=\Cal L_W$  for any $W\in \Cal
 F$.  As proved in \cite{PS, \S7.7}, the line bundle $\Cal L$
 is an algebraic line bundle on $\Cal F$ (with respect to
 the ind-variety structure on $\Cal F$  as in \S C.7).
 It is easy to see that $\Cal L_{\vert_{\Cal F_1}}$  is
 the restriction of the basic (negative ample) line
 bundle on $Gr(nN,2nN)$ under the identification $\Cal
 F_1 \overset\sim\to{\rightarrow}
 Gr(nN,2nN)^{1+\overline{t}_1}$ (cf. \S C.7).  Let $\Cal
 L_o$ be the pull-back of the line bundle $\Cal L$ to
 $X^o$ via the isomorphism $\beta :X^o
 \overset\sim\to{\rightarrow} \Cal F$.  In view of
 Proposition (C.13), it is easy to see that the dual line
 bundle $\Cal L^\ast_o$ is isomorphic with the line
 bundle ${\frak L}(\chi_0)$.

 Now we define an action $\alpha$ of the group $\Cal H/\Cal N$
 on $\Cal L$ as follows:  For $(g,E)\in \Cal H$, define
 $$
 \alpha(g,E) [z,\theta]_W =
 [z,g\theta E^{-1}]_{gW} ,
 $$
 where for $z\in \Bbb C$ and $\theta\in S_W$,
 $[z,\theta]_W$ denotes the equivalence class of
 $z\theta$.  This action factors through an action of
 $\Cal H/\Cal N$ and moreover for any fixed $(g,E)\in
 \Cal H$, $\alpha(g,E)$ is an algebraic automorphism of
 the line bundle $\Cal L$ (and hence of $\Cal L_o$)
 inducing the map $L_g$ on the base (cf. \S C.6).  Using this, the group $\Cal
 H/\Cal N$ can canonically be identified with the Mumford
 group Aut$( {\frak L}_o)$ defined in \S (C.6). In
 particular, the central extension $\Cal H/\Cal N$ is
 isomorphic with $\tilde{\Cal G}$.

 Finally we construct a splitting of $\Cal H/\Cal N$ over
 $\Gamma$ as follows:

Choose an element $g_o\in \Cal G^o$ such that the
 associated rank-$N$ vector bundle $\tilde{\varphi} (g_o)$ on $C$ twisted by
$\Cal O((g-1)p),
 E:= \tilde{\varphi} (g_o)((g-1)p)
 $
 (where $g$ is the genus of the curve $C$) has all
 its cohomology $0$. Then considering  the local cohomology sequence
 (for the curve $C$ with support in $p$) with coefficients in  the
 vector bundle $E$, we deduce that
 $$
 V((t)) = L_o \oplus t^{1-g}
 g_o^{-1}(V\otimes_{\Bbb C}{\Bbb C}[C^{\ast}])~,
 \tag"(3)" $$
 where $V\otimes_{\Bbb C}{\Bbb C}[C^\ast]$ is
 identified as a subspace of $V((t))$ by choosing a
 parameter $t$ around $p \in C$.

 We first construct the splitting of $\Cal H/\Cal N$
over
 $\Gamma_{g_o}:=g_o^{-1}\Gamma g_o$.
 Define the group homomorphism $\sigma_{g_o}:
 \Gamma_{g_o}\to\Cal H$ by $\sigma_{
 g_o}(\gamma)=(\gamma,\gamma^{+'}) $, where $\gamma
 =\pmatrix\gamma^{+'} &0\\ \ast &\ast\endpmatrix$ with
 respect to the decomposition (3).
 (Observe that $\Gamma_{g_o}$ keeps the second
 factor stable and hence $\gamma^{+'}\in GL(L_o)$.)  The
 group homomorphism $\overline{\sigma}_{g_o}: \Gamma
_{g_o} \to \Cal H/\Cal N$ (where
 $\overline{\sigma}_{g_o}$ is the map $\sigma_{
 g_o}$ followed by the canonical map $\Cal H\to \Cal
 H/\Cal N$) splits the central extension (2) over $\Gamma
_{g_o}$.  Now take any preimage $\overline{
 g}_o$ of $g_o$ in $\Cal H/\Cal N$, and define the
 splitting $\overline{\sigma}: \Gamma \to \Cal H/\Cal N$
 ($\gamma \mapsto\overline{g}_o \overline{\sigma}_{g_o}
 (g_o^{-1}\gamma g_o)\overline{g}_o^{-1}$).

 Since $\Cal H/\Cal N$ acts on the line bundle $\Cal
 L_o$, so is $\Gamma$ (via the homomorphism
 $\overline{\sigma }$).  It can be easily seen that the
 action $\Gamma \times \Cal L_o\to \Cal L_o$ is a morphism of ind-varieties.
Moreover, let $\overline{\sigma}' : \Gamma \to \Cal H/\Cal N $ be another
splitting of $\Gamma$ such that the induced action  $\Gamma \times \Cal L_o\to
\Cal L_o$ is again  a morphism of ind-varieties. Then we claim that
$\overline{\sigma}' = \overline{\sigma}$ :
There is a group homomorphism $\alpha : \Gamma \to \Bbb C^\ast$ such that
(cf.(2)) $\overline{\sigma}'=\alpha \overline{\sigma}$. Further $\alpha$
is a morphism of ind-varieties (since the action of $\Gamma$ on $\Cal L_o$ in
both the cases is regular). But then $\alpha$ is identically $1$ (cf. Remark
6.8(c)). This proves the uniqueness of such a splitting.

Since the line bundle $\Cal L_o$ is isomorphic with the homogeneous line bundle
$\frak L(-\chi_0)$, it is easy to see that the group $\Gamma$ acts morphically
on the representation $L(\Bbb C, 1)$ and hence on any $L(\Bbb C, d)$ (for $d
\geq 0$, where $L(\Bbb C, d)$ is the irreducible representation of the affine
Lie algebra $\tilde{sl_N}$ with central charge $d$, cf. \S A.2).

\flushpar \subhead (2.7) Splitting of the central
 extension over $\Gamma$ (general case) \endsubhead
We now come to the case of general $G$ as in \S1.1.
Take a finite dimensional representation $V$ of $G$
such that the group homomorphism $\varphi : G \to \text{SL}(V)$ has finite
kernel, and consider the
induced map at the Lie algebra level $d\varphi: \frak g \to sl(V)$. We denote
the Lie algebra $sl(V)$ by $\frak g^o$. The Lie algebra homomorphism
 $d\varphi$ induces a Lie algebra homomorphism $\tilde{\varphi} : \tilde{\frak
g} \to \tilde{\frak g}^o$ defined by (cf. \S A.1)
$$ X\otimes p \mapsto (d\varphi (X))\otimes p~,~
{}~\text{and}~~ K\mapsto m_VK^o~ ,$$
for $X \in \frak g$ and $p\in \Bbb C[t^{\pm 1}]$; where $K$ (resp. $K^o$) is
the canonical central element of  $\tilde{\frak g}$ (resp. $\tilde{\frak
g}^o$), and $m_V$ is the Dynkin index of the representation $V$ (cf. \S5.1).

To distinguish the objects corresponding to SL$(V)$ from that of $G$, we denote
the former by a superscript $o$.
Let us consider the irreducible representation
$L^o(\Bbb C, 1)$ for the Lie algebra $\tilde{\frak g}^o$ with central charge
$1$ and restrict it to the Lie algebra
$\tilde{\frak g}$ via the homomorphism $\tilde{\varphi}$. It can be seen that
the
$\tilde{\frak g}$-submodule of $L^o(\Bbb C, 1)$ generated by the highest weight
vector $v_o$ is isomorphic with $L(\Bbb C, m_V)$.

The representation $\varphi$ also gives
 rise to a morphism of the corresponding affine flag
varieties $\hat{\varphi}: X \to X^o$, and a morphism of ind-groups $ \Gamma \to
\Gamma^o$. It is easy to see that the basic homogeneous line bundle $\frak
L^o(\chi_0) $ on $X^o$ pulls-back to the line bundle $\frak
L(m_V\chi_0) $ on $X$. In particular, the group $\Gamma$ acts morphically on
the line bundles $\frak
L(dm_V\chi_0) $ (for any $d \in \Bbb Z$) and hence $\Gamma$ also acts
morphically on the representation space $L(\Bbb C, dm_V)$.

\vskip1ex
Finally we have  the following proposition, asserting that
 $X= \Cal G/\Cal P$ supports an algebraic family.

\proclaim{(2.8) Proposition} (a) There is an algebraic $G$-bundle
 $\Cal U \to C\times X$  such that
 for any $x\in X$ the $G$-bundle $\Cal U_x :=
\Cal U_{\mid C\times x}$ is isomorphic with $\varphi(x)$ (where
 $\varphi$ is the map of \S1.4). Moreover,  the bundle
$\Cal U_{\mid C^{\ast} \times X}$ comes equipped with a
trivialization $ \alpha : \epsilon\, \tilde{\rightarrow}  \,\Cal U_{\mid
C^{\ast} \times X}$, where $\epsilon$ is the trivial $G$-bundle on
$C^{\ast} \times X$.

(b) Let $\Cal E \to C\times Y$ be an algebraic  family of $G$-bundles
(parametrized
by an algebraic variety $Y$), such that $\Cal E$ is trivial over
$C^\ast \times Y$ and also over    (Spec $\hat {\Cal O}_p) \times Y$.
Then, if we choose a trivialization $\beta : \epsilon' \, \tilde{\rightarrow}
\,\,\Cal E_{\mid  C^\ast \times Y}$, we get a Schubert variety $X_{\frak w}$
and a unique morphism $f : Y \to X_{\frak w}$ together with
a $G$-bundle morphism $\hat{f} : \Cal E \to \Cal U_{\mid C\times X_{\frak w}}$
inducing the map I$\times f$ at the base such that $\hat{f} \circ \beta =
\alpha \circ \theta$, where $\epsilon'$ is the trivial bundle on $C^\ast \times
Y$ and $\theta$ is the canonical $G$-bundle morphism  $\epsilon' \to
\epsilon$ inducing the map I$\times f$ at the base.
\endproclaim
\demo{Proof} Let $R$ be a $\Bbb C$-algebra and let  $Y$ := Spec $R$ be the
corresponding scheme.  Suppose $E \to C \times Y$ is a $G$-bundle with
trivialisations $\alpha$
of $E$ over $C^\ast \times Y$ and $\beta$ of $E$ over    (Spec $\hat {\Cal
O}_p) \times Y$. Note that
the fiber product $(C^\ast \times Y) \times_{C\times Y} (\text{Spec}\,
\hat{\Cal O}_p \times Y)$  is canonically isomorphic with ~~~~~~
(Spec $\hat{\text{{\sl k}}}_p) \times Y$.  Therefore the trivialisations
$\alpha$ and $\beta$ give rise to
an element $\alpha \beta^{-1} \in G(\hat{\text{{\sl k}}}_p \otimes R)$.
Conversely, given
an element $g \in G(\hat{\text{{\sl k}}}_p \otimes R)$, we can construct the
family $E \to C \times \text{Spec} ~R$ by taking the trivial bundles on
$C^\ast \times Y$ and (Spec $\hat{\Cal O}_p) \times Y$ and glueing
them via the element $g$. Moreover,
 if $g_1$ and $g_2$ are two elements of  $G(\hat{\text{{\sl k}}}_p \otimes R)$
such that $g_2 = g_1 h$ with $h \in G(\hat{\Cal O}_p \otimes R$), then $h$
induces a canonical
isomorphism of the bundles corresponding to $g_1$ and $g_2$.  All these
assertions
are easily verified.

Choose a local parameter $t$ around $p\in C$. Let
ev$_\infty : G(\Bbb C[t^{-1}]) \to G$ be the group
homomorphism induced from the algebra homomorphism
$\Bbb C[t^{-1}] \to \Bbb C $ taking $t^{-1} \mapsto 0$ , and let $N^- :=$ ker
(ev$_\infty$). Then the image $U^-$ of $N^-$ in $X$ under the map $i: N^- \to X
$ taking $g\mapsto g\Cal P $ is an open subset of $X$.
 To construct a family of $G-$bundles on $X,$
 we first construct a family on the open sets $ w U^-  \subset X $, for $ w \in
\text{Mor}_1 (\Bbb C^*, T)$ as follows (cf. proof of Lemma C.10 for the
notation $\text{Mor}_1 (\Bbb C^*, T)$):

 From the discussion in the first paragraph, it suffices to construct an
element ${\theta}_w \in G( \hat{\text{{\sl k}}}_p \otimes  k[w U^-])$ such
 that for every $w x \in w U^-$, the element ${\theta}_w$
 evaluated at $w x$ (i.e. the image of ${\theta}_w$ under the
 evaluation map  $G (\hat{\text{{\sl k}}}_p \otimes   k[w U^-])
 \rightarrow G (\hat{\text{{\sl k}}}_p)$ at $wx)$ satisfies
 ${\theta}_w (w x)= w i^{-1}(x) ~\text{mod} ~{\Cal P}$.   But,  by definition,
$N^- \subset  G (\Bbb C [t^{-1}])$ and
 hence we get a tautological map
 ${\theta}: (\Bbb P^1 (\Bbb C) \backslash 0) \times N^- \rightarrow G$.  It is
  easy to see that ${\theta}$ is a morphism under the
 ind-variety structure on $N^-$. (Observe that $U^-$ being an open subset of
$X_{\text{lat}}$ has an ind-variety structure and hence $N^-$ acquires an
ind-variety structure via the bijection $i$.) Think of $  \Bbb C^* = \Bbb P^1
(\Bbb C) \backslash \{
 0, \infty\} $ and  define
 $\overline{{\theta}}_w: \Bbb P^1 (\Bbb C) \backslash \{ 0, \infty\}
 \times w U^- \rightarrow G$ by
 $\overline{{\theta}}_w (z, w i(g))=w(z) {\theta} (z,g),$
 for $z \in \Bbb P^1 (\Bbb C) \backslash \{ 0, \infty\}$ and  $g \in
 N^-$.  The morphism $\overline{{\theta}}_w $ of course gives
 rise to an element
 ${\theta}_w \in G(\hat{\text{{\sl k}}}_p \otimes  k[w
 U^-]),$  and hence a $G$-bundle on $C \times wU^-$.

 To prove that the bundles on $C \times w U^-$ got from the
 elements ${\theta}_w$ patch together to give a bundle on $C
 \times X$, it suffices to show that the map
 $$\overline{{\theta}}_{v}^{-1} \overline{{\theta}}_w: \Bbb P^1
 (\Bbb C) \backslash \{0, \infty\} \times (w U^- \cap  v U^-)
 \rightarrow G$$
 extends to a morphism (again denoted by)
 $\overline{{\theta}}_{v}^{-1} \overline{{\theta}}_w: \Bbb P^1
 (\Bbb C) \backslash \{ \infty\} \times (w U^- \cap  v U^-)
 \rightarrow G:$
 But for any fixed $x \in w U^- \cap  v U^-,$ the map
 $\overline{{\theta}}_{v}^{-1}\overline{{\theta}}_{w}: \Bbb P^1
 (\Bbb C) \backslash \{ 0, \infty \} \times x \rightarrow G$  in fact
 is an element of ${\Cal P}=G  (\widehat{\Cal
 O}_p)$, i.e.,  $\overline{{\theta}}_{v}^{-1} \overline{{\theta
 }}_{w}$ does not have a pole at $0$, for any fixed $x \in w U^-
 \cap  v U^-$.  From this it is easy to see that $\overline{{\theta
 }}_{v}^{-1} \overline{{\theta}}_{w}$ extends to a morphism
 $\Bbb P^1 (\Bbb C) \backslash \{\infty\} \times (w U^- \cap  v U^-)
 \rightarrow G$.  Clearly the maps $\overline{{\theta}}_{v}^{-1}
 \overline{{\theta}}_{w}$ satisfy the `cocycle condition' and hence
 we get a $G-$bundle on the whole of $C \times X$.

To prove the (b) part, let us choose a trivialization $\tau$ of the bundle
$\Cal E$ restricted to $(\text{Spec}~ \hat {\Cal O}_p) \times Y$. As above,
this (together with the trivialization $\beta$) gives rise to a map $f_\tau :
Y\to \Cal G $ and hence  a map
$f : Y\to X.$ (It is easy to see that the map $f$ does not depend upon the
choice of the trivialization $\tau$.) We claim that there exists
a large enough $X_{v}$ such that Im$ f \subset X_{\frak w}$ and
moreover $f : Y \to X_{\frak w}$ is a morphism:

For both of these assertions, we can assume that $Y$ is an affine variety
$Y= \text{Spec}~ R$, for some $\Bbb C$-algebra $R$. Then the map $f_\tau$ can
be thought of as an element (again denoted by) $f_\tau \in G(\hat{\text{{\sl
k}}}_p \otimes R)$. Choose an imbedding $G \hookrightarrow SL_N$.
  Then we can
write $f_\tau = (f_\tau^{i,j})_{ 1\leq i,j \leq N}$, with
$f_\tau^{i,j} \in \hat{\text{{\sl k}}}_p \otimes R$. In particular,
 there exists a large enough $l \geq 0$ such that (for any
 $1\leq i,j \leq N$) $f_\tau^{i,j} \in t^{-l} \Bbb C[[t]] \otimes R$. From this
(together with Lemma C.10) one can see that
 Im $f$ is contained in a Schubert variety $X_{\frak w}$. Now
the assertion that $f : Y\to X_{\frak w}$ is a morphism follows from
the description of the map $f^\tau$ as an element of  $G(\hat{\text{{\sl k}}}_p
\otimes R)$ together with the  description of the variety structure
$X_{\text{lat}}$
on $X$. The remaining assertions of (b) are easy to verify, thereby completing
the proof of  (b). \qed
\enddemo

\vskip6ex
\subheading {3. Preliminaries on moduli space of
 $G$-bundles and the determinant bundle}
\vskip4ex

  {\it Throughout this section, we allow $G$  to be  a connected
 reductive group over an algebraically closed field
 $k$  of char. 0. }

 We recall some basic concepts and results on
 semistable $G$-bundles on $C$.  The references are
 \cite{NS}, \cite{R$_1$}, \cite{R$_2$}, and \cite{RR}.
 Recall the definition of $G$-bundles and reduction of
 structure group from \S 1.2.
\vskip1ex

 \flushpar {\bf (3.1)} {\it  Definition}. Let $E\to C$ be a
 $G$-bundle.  Then $E$ is said to be {\it semistable}
 (resp. {\it stable}), if for any reduction of
 structure group $E_P$ to any parabolic subgroup
 $P\subset G$ and any non-trivial character $\chi  :P \to {\Bbb G}_m$
 which is dominant with respect to some Borel subgroup
 contained in $P$, the degree of the associated line
 bundle $E_P(\chi)$ is $\leq 0$ (resp. $<0$).  (Note
 that, by definition, a dominant character is taken to
 be trivial on the connected component of the centre
 of $G$.)

 \flushpar {\bf (3.2)}
{\it  Remark}.  When $G=GL_n$, this
 definition coincides with the usual definition of
 semistability (resp. stability) due to Mumford (cf.
 \cite{NS}) vz. a vector bundle $V\to C$ is
 semistable (resp. stable) if for every subbundle
 $W\subsetneqq V$, we have $\mu(W)\leq \mu(V)$ (resp.
 $\mu(W) < \mu(V)$), where $\mu(V):= \text{deg}\,V/\text{rank}\,V$.

 Let $V\to C$ be a semistable vector bundle.  Then
 there exists a filtration by subbundles
 $$
 V_0=0\,\subsetneqq V_1 \subsetneqq
 V_2 \subsetneqq \dots \subsetneqq V ,
 $$
 such that $\mu(V_i)=\mu(V)$ and $V_i/V_{i-1}$ are stable.
  Though such a filtration in
 general is not unique, the associated graded
 $$
 \text{gr}\,V := \underset{i\geq 1}\to\oplus\, V_i/V_{i-1}
 $$
 is uniquely determined by $V$ (upto an isomorphism).

 We will now describe the corresponding notion of
 gr$E$ for a semistable $G$-bundle $E$.

 \flushpar {\bf (3.3)} {\it  Definition.} A reduction of structure
 group of a $G$-bundle $E\to C$ to a parabolic
 subgroup $P$ is called {\it admissible} if for any
 character of $P$, which is trivial on the connected
 component of the centre of $G$, the associated line
 bundle of the reduced $P$-bundle has degree 0.

 It is easy to see that if $E_P$ is an admissible
 reduction of structure group to a parabolic subgroup
 $P$, then $E$ is semistable if and only if the
 $P/U$-bundle $E_P(P/U)$ is semistable, where $U$ is
 the unipotent radical of $P$.  Moreover, a semistable
 $G$-bundle $E$ admits an admissible reduction to some
 parabolic subgroup $P$ such that $E_P(P/U)$ is, in
 fact, a stable $P/U$-bundle.  Let $M$ be a Levi
 component of $P$.  Then $M \approx P/U$ (as algebraic
 groups) and thus we get a stable $M$-bundle $E_P(M)$.
 Extend the structure group of this $M$-bundle to $G$
 to get a semistable $G$-bundle denoted by gr$(E)$.
 Then gr$(E)$ is uniquely determined by $E$ (up to an
 isomorphism).

 Two semistable $G$-bundles $E_1$ and $E_2$ are said
 to be $S$-{\it equivalent} if gr$(E_1)\approx$ gr$(E_2)$.
 We call a semistable $G$-bundle $E$ {\it quasistable}
 if $E\approx \text{gr}(E)$.  (It can be seen that a
 semistable vector bundle is quasistable if and only
 if it is a direct sum of stable vector bundles with the
 same $\mu$.)

  Two $G$-bundles $E_1$ and
 $E_2$ on $C$ are said to be of the {\it same
 topological type} if they are isomorphic as
 $G$-bundles in the topological category.  The
 topological types of all the algebraic $G$-bundles
 are bijectively parametrized by $\pi_1(G)$ (cf.
 \cite{R$_2$, \S 5}).

 \proclaim{(3.4)  Theorem}  The set $\frak M$ of
 $S$-equivalence classes of all the semistable
 $G$-bundles of a fixed topological type admits the
 structure of a normal, irreducible, projective
 variety over $k$, making it into a coarse moduli.

 In particular, for any algebraic family $\Cal E \to C\times Y$ of
 semistable $G$-bundles of the same topological type
 (parametrized by a variety $Y$), the set map $\beta:Y\to
 \frak M$, which takes $y\in Y$ to the $S$-equivalence
 class of $\Cal E_y$ in $\frak M$ is a morphism.
 \endproclaim

 The details can be found in \cite{NS}, \cite{R$_1$},
 \cite{R$_2$}, ....

 \proclaim{(3.5) Lemma}
Let $H$ be a connected affine algebraic group and $C$ a smooth
projective curve over $k$. Then any principal $H$-bundle on $C$ is
locally trivial in the Zariski topology.
\endproclaim

\demo{Proof} Let $E$ be a principal $H$-bundle on $C$ and $U$
the unipotent radical of $H.$  Since the group $M = H/U$ is
connected and reductive, the $M$-bundle $E(M),$ obtained from
$E$ by extension of structure group to $M,$ is locally trivial
in the Zariski topology  \cite{$\text{R}_3$, Proposition 4.3}.

Let $W$ be a non-empty affine open subset of $C$ such that the
restriction of $E(M)$ to $W$ is trivial. We  show that $E_
{\mid W}$ is trivial (which  will prove the lemma): Observe that a
trivialisation of $E(M)$ on $W$ gives a reduction of the
structure group $H$ of $E_{\mid W}$ to the subgroup $U.$ So, it suffices
to show that any (principal) $U$-bundle on $W$ is trivial:

We may
assume $U \neq e.$ Then there exists a (finite) filtration of $U$
by closed normal subgroups such that the successive quotients
are isomorphic to the additive group $G_{a}.$  Now the assertion follows since
any
principal $G_{a}$-bundle on $W$ is trivial, $W$ being
affine (see \cite{$\text{Se}_1$, \S5.1}).
\enddemo

Let $P$ be a parabolic subgroup of $G$ and $P = MU$ a Levi
decomposition,  where $U$ is the unipotent radical of $P$
and $M$ a Levi component. The next proposition will be used in \S6 in the case
of an
admissible reduction of a  semistable bundle $E.$

\proclaim{(3.6) Proposition} Let $G$ be a connected semisimple algebraic group.
Let $E$ be a $G$-bundle on $C$ and
$E_{P}$ a reduction of the structure group  of $E$ to $P.$
Denote by gr $ (E_{P})$ the $G$-bundle on $C$ obtained from the
$P$-bundle $E_{P}$ by extension of the structure group
via the composite homomorphism
$$P \rightarrow P/U \approx M \hookrightarrow G.$$

Then there exists a $G$-bundle ${\Cal E}$ on $C\times \Bbb A^1,$ where
${\Bbb A}^1$ is the affine line, such that we have
\roster
\item "(a)" ${\Cal E}_{\mid C\times (\Bbb A^1 \backslash 0) } \approx
p^{\ast}_{C} (E), \ {\Cal E}_{\mid C \times \{0\} } \approx
 \text{gr}\ (E_{P})$
and
\item "(b)" ${\Cal E}_{\mid C^\ast \times \Bbb A^1}$ is trivial
and also the pull-back of
 ${\Cal E}$ to  (Spec $\hat{\Cal O}_p) \times \Bbb A^1$
 is trivial,
\endroster
where $p_C$ is the projection on the $C$-factor.
\endproclaim

\demo{Proof} By  \cite{$\text{R}_1$, Lemma 2.5.12}, there exists
 a one-parameter
group $\lambda : G_{m} (:= \Bbb A^1 \backslash 0) \to  M$ ,  such that
the  regular map
$$G_{m} \times P \rightarrow P,\,\,  \text{given\, by}\, (t,p) \mapsto
 \lambda (t) p
\lambda (t)^{-1},\,\, \text{for} \,\, t \in G_{m},\, p \in P,$$
 extends to a regular map  $\phi : \Bbb A^1 \times P \rightarrow P$
satisfying  $\phi (0, mu) = m$ , for $  m
\in M, u \in U.$  By Lemma (3.5), the $P$-bundle $E_{P}$
is locally trivial in the Zariski topology. Let $\{
U_{i}\}$ be an affine open covering of $C$ in which the bundle
$E_{P}$ is given by the transition functions $p_{ij} : U_{i} \cap
U_{j} \rightarrow P.$ Let ${\Cal F}$ be the (Zariski locally
trivial) $P$-bundle on $C \times \Bbb A^1$ defined by the covering
$\{U_i \times \Bbb A^1 \}$ and the transition functions
$$h_{ij} :  (U_{i} \cap U_{j}) \times \Bbb A^1 \rightarrow P \,,
$$
where $h_{ij}(z,t) = \phi(t,p_{ij}(z)),$ for $ t \in \Bbb A^1, z \in U_{i} \cap
U_{j}.$ Now let  ${\Cal E}$ be the $G$-bundle obtained from the $P$-bundle
 ${\Cal
F}$ by extension of the structure group to $G$. Then clearly  ${\Cal E}$
 satisfies condition
(a).

We next   show that for any non-empty affine  open
subset $W$ of $C,$ the restriction of ${\Cal E}$ to $W \times \Bbb A^1$ is
trivial (this will, in particular, imply that condition b is
satisfied): Note that, by our construction,  there exists a finite
open covering $W_{i}$ of $W$ such that ${\Cal E}_{ \mid W_i \times
\Bbb A^1}$ is trivial, for every $i.$ Now by an analogue of a result
of Quillen (cf. \cite{Ra, Theorem 2})  ${\Cal E}_{ \mid W \times \Bbb A^1}$
is the pull-back of a $G$-bundle on $W.$ But by Proposition (1.3),
any  $G$-bundle on $W$ is trivial.
\enddemo

 \vskip1ex

 \flushpar {\bf (3.7)} {\it   Determinant bundle and $\Theta$-bundle}.
  We now briefly recall a few definitions and facts on the
determinant bundles and $\Theta$-bundles associated to families of bundles on
$C$.  We follow \cite{DN}, \cite{NRa}.

Let $\Cal V \to C\times Y$ be a vector bundle.  Then there exists a complex
of vector bundles  $\Cal V_i $ on $Y$ ( with $\Cal V_i  = 0$, for all $i \geq
2$):
$$
\Cal V_0\to \Cal V_1\to 0 \to 0 \to \dots,
$$
 such that for any base change $f:Z\to Y$, the
$i^{\text{th}}$ direct image (under the projection $C\times Z \to Z$) of the
pull-back$(\text{id}\times f)^\ast \Cal V$ on
$Z$ is given by the $i^{\text{th}}$ cohomology of the pull-back of the above
complex to $Z$.  We define the {\it determinant line bundle} Det $\Cal V$ on
$Y$ to be the
 product $\overset\text{top}\to{\wedge}(\Cal V_{1})\otimes
(\overset\text{top}\to{\wedge}(\Cal V_{0})^\ast)$. ( Notice that our
Det $\Cal V$ is dual to the determinant line bundle as defined, e.g., in
\cite{L, Chapter 6,\S 1}.)

 The above base change property gives rise to the base change property for
Det $\Cal V$, i.e.,  if $f: Z \to Y$ is a morphism then
Det$ ((\text{id}\times f)^\ast\Cal V)= f^\ast$(Det $\Cal V)$.

Let  $\frak L$ be  a line bundle on $Y$, and let $p_2: C\times Y \to Y$
 be the projection on the second factor. Then for the family $\Cal V \otimes
p_2^\ast
\frak L\to C\times Y$, we have Det $(\Cal V\otimes p_2^\ast \frak L)$= (Det
$\Cal V)\otimes
\frak L^{-\Cal X(\Cal V)}$, where $\Cal X(\Cal V)
 :=h^0(\Cal V_t)-h^1(\Cal V_t)$ is the Euler characteristic and $\Cal V_t :=
\Cal V \vert_{C\times t}$. (Observe that $h^0(\Cal V_t)-h^1(\Cal V_t)$ remains
constant on any connected component of $Y$.)

We now define the $\Theta$-{\it bundle} $\Theta(\Cal V)$ of a family
of rank $r$ and degree $0$ bundles $\Cal V \to
C\times Y$ to be the modified determinant bundle given by (Det $\Cal V)\otimes$
(det$(\Cal V_p))^{\Cal X(\Cal V)/r}$, where $\Cal V_p$ is the bundle $\Cal V
\vert_{p\times Y}$ on $Y$, and det $\Cal V_p$ is its usual determinant line
bundle.  It follows then that $\Theta(\Cal V)=\Theta(\Cal V \otimes
p^\ast_2\frak L)$, for any line bundle $\frak L$ on $Y$. Moreover $\Theta(\Cal
V)$ also has the functorial property
$\Theta((\text{id}\times f)^\ast\Cal V) = f^\ast (\Theta(\Cal V))$.

If $\Cal E\to C\times Y$ is a family of $G$-bundles ( where $G$ is semisimple
and connected) and  $V$ is
a $G$-module, then Det $(\Cal E(V))$ and $\Theta(\Cal E(V))$ are defined to be
the corresponding  line bundles of the associated family of vector bundles,
via the representation $V$ of $G$.

For the  family $\Cal U \to C\times X$ (cf. Proposition 2.8),
 the line bundles $\Theta(\Cal U(V))$ and $\text{Det} (\Cal U(V))$ coincide,
since
$\Cal U_{\vert_{p\times \Cal G/\Cal P}}$ is trivial.

 It is known (\cite{DN}, \cite{NRa} ) that there exists a line bundle
$\Theta$ on the moduli space $\frak M_o$ of rank $r$ and degree $0$
(semistable) bundles, such that for any family $\Cal V$ of rank $r$ and degree
$0$ semistable bundles parametrized by $Y$ we have $f^{\ast}(\Theta ) \simeq
\Theta (\Cal V)$, where $f: Y \to \frak M_o$ is the morphism given by the
coarse moduli property of $\frak M_o$ (cf. Theorem 3.4).

Let $V$ be a  representation of $G$ of dimension $r$ ($G$ semisimple and
connected).
 Then for any semistable $G$-bundle on $C$, the associated vector bundle
(via the representation $V$)  is  semistable (cf. \cite{RR, Theorem 3.18}).
Thus, given a family of semistable $G$-bundles on $C$ parametrized by
$Y$, we have a canonical morphism (induced from the representation $V$) $Y \to
\frak M_o$ (where $\frak M_o$
as above is the  moduli space of semistable bundles of rank $r$ and degree $0$
). Let $\frak M$ be the moduli space of semistable $G$-bundles. By the coarse
moduli property of $\frak M$, we see that we have a canonical morphism
$\phi_V : \frak M \to \frak M_o$. We define the {\it theta bundle $\Theta (V)$
on $\frak M$ associated to $V$} to be the pull-back of the line bundle
$\Theta$ on $\frak M_o$ via the morphism  $\phi_V$.
 It can be easily seen that for any family $\Cal V \to C \times Y$
of semistable $G$-bundles, $f^{\ast} (\Theta (V)) \simeq \Theta (\Cal V(V))$,
where $f : Y \to \frak M$ is the morphism (induced from the family $\Cal V$)
given by the coarse moduli property of $\frak M$.
\vskip5ex
 \subheading{ 4.  A result on algebraic descent}
\vskip2ex
We prove the following technical result, which will crucially be used in the
note.  Even though we believe that it should be known, we did not find a
precise
reference.
\proclaim{(4.1) Proposition}  Let $f:X\to Y$ be a surjective morphism between
 irreducible algebraic varieties $X$ and $Y$ over an algebraically closed field
$k$ of char
$0$. Assume that $Y$ is normal and let $\Cal E
\to Y$ be an algebraic vector bundle on $Y$.

Then any set theoretic section $\sigma$ of the vector bundle $\Cal E$ is
regular
if and only if the induced section $f^\ast(\sigma)$ of the induced bundle
$f^\ast(\Cal E)$ is regular.
\endproclaim
\demo{Proof}  The `only if' part is of course trivially true.  So we come to
the
`if' part.

Since the question is local (in $Y$), we can assume that $Y$ is affine and
moreover the vector bundle $\Cal E$ is trivial, i.e., it suffices to show that
any (set theoretic) map $\sigma : Y \to k$ is regular, provided $\bar{\sigma}:=
\sigma \circ f: X\to k$ is regular (under the assumption that $Y =$ Spec $R$ is
irreducible normal and affine):

Since the map $f$ is surjective (in particular dominant), the ring $R$ is
canonically embedded in $\Gamma (X):= H^0(X,\Cal O_X)$.  Let $R[\bar{\sigma}]$
denote the subring of $\Gamma (X)$ generated by $R$ and $\bar{\sigma}\in \Gamma
(X)$.  Then $R[\bar{\sigma}]$ is a (finitely generated) domain (as $X$ is
irreducible by assumption), and we get a dominant morphism $\hat{f}: Z\to$ Spec
$R$, where $Z:=$ Spec $(R[\bar{\sigma}])$.  Consider the commutative diagram:
$$
\gather
X\\
\theta\swarrow\qquad \searrow f\\
Z \underset\hat{f}\to\longrightarrow Y
\endgather
$$
where $\theta$ is the dominant morphism induced from the inclusion
$R[\bar{\sigma}]\hookrightarrow \Gamma (X)$.  In particular, Im $\theta$
contains
a non-empty Zariski open subset $U$ of $Z$.  Let $x_1,x_2\in X$ be closed
points such that $f(x_1)= f(x_2)$.  Then $r(x_1)= r(x_2)$, for all $r\in R$ and
also $\bar{\sigma}(x_1)= \bar{\sigma}(x_2)$.  This forces $\theta(x_1)
=\theta(x_2)$, in particular, $\hat{f}_{\vert U }$ is injective on the closed
points of $U$.

Since $\hat{f}$ is dominant, by cutting down $U$ if necessary, we can assume
that $\hat{f}_{\vert U }: U \to V$ is a bijection, for some open subset
$V\subset Y$.  Now since $Y$ is (by assumption) normal and $Z$ is irreducible,
by Zariski's main theorem (cf. \cite{Mum, Page 288, I. Original form}),
$\hat{f}_{\vert U }:U \to V$
is an isomorphism, and hence $\sigma$ is regular on $V$.

 Assume, if possible, that $\sigma_{\vert_V}$ does not
extend to a regular function on the whole of $Y$.  Then, by \cite{B, Lemma
18.3, Chapter AG},
there exists a point $y_0\in Y$ and a regular function $h$ on a Zariski
neighborhood $W$ of $y_0$ such that $h(y_0)=0$ and $h\sigma=1$ on $W\cap V$.
But then $\bar{h}\bar{\sigma}=1$ on $f^{-1}(W\cap V)$ (where $\bar{h}:= h\circ
f$) and hence, $\bar{\sigma}$ being regular on the whole of $X$,
$\bar{h}\bar{\sigma}=1$ on $f^{-1}(W)$.  Taking $\bar{y_0}\in f^{-1}(y_0)$ ($f$
is,
by assumption, surjective), we get
$\bar{h}(\bar{y_0})\bar{\sigma}(\bar{y_0})=0$.  This contradiction shows that
$\sigma_{\vert_V}$ does extend to some regular function (say $\sigma'$) on the
whole of $Y$.  Hence $\bar{\sigma}= \bar{\sigma}'$, in particular, by the
surjectivity of $f$, $\sigma=\sigma'$.  This proves the proposition.
\qed \enddemo

\vskip4ex

 \subheading{ 5. Identification of the determinant
 bundle}
\vskip5ex

{\bf (5.1) } Recall from \S 2.8 that $\Cal G/\Cal P$ is a parameter space for
an algebraic family $\Cal U$ of $G$-bundles on $C$. Let us fix a (finite
 dimensional) representation $V$ of $G$. In particular, we can talk of the
determinant line  bundle  Det$(\Cal U(V))$ (cf. \S 3.7). Also recall
the definition of the fundamental homogeneous line bundle $\frak L(\chi_o)$
on $\Cal G/\Cal P$ from \S (C.6). Our aim in this section is to determine
the line bundle Det$(\Cal U(V))$ in terms of  $\frak L(\chi_o)$. We begin
with the following preparation.

Let $\theta$ be the highest root of $\frak g$. Define the following Lie
subalgebra $sl_2(\theta)$
  of the Lie algebra $\frak g$ of $G$ :
$$ sl_2({\theta}) := \frak g_{-{\theta}}
 \oplus {\Bbb C}
{\theta}^\vee \oplus\frak g_{\theta} , \tag"(1)" $$ where
 $\frak g_{\theta} $ is the ${\theta}$-th root
 space, and ${\theta}^\vee$ is the corresponding coroot.
 Clearly $sl_2({\theta}) \approx sl_2$ as Lie
 algebras. Decompose $$ V=\oplus_{i}  V_i ,
 \tag"(2)" $$ as a direct sum of irreducible
 $sl_2 ({\theta})$- modules $V_i$ of dim $m_i$. Now
 we define $$ m_V=\sum_i \left(\matrix
 m_i+1\\3\endmatrix\right)\ , ~~~ \text{ where we set}\ \pmatrix
 2\\3\endpmatrix=0.\tag"(3)"$$

  Let $\frak{g}_1$ and
 $\frak{g}_2$ be two (finite dimensional) complex simple
 Lie algebras and $\varphi : \frak{g}_1 \to \frak{g}_2$
 be a Lie algebra  homomorphism.  There exists a unique
 number $m_\varphi \in \Bbb C$, called the {\it Dynkin index} (cf. \cite{D,
\S2})
 of the homomorphism $\varphi$, satisfying
 $$
 \langle \varphi (x),\varphi (y) \rangle
 = m_\varphi \langle x,y \rangle, \text{ for all }x,y \in
 \frak{g}_1,
 $$
 where $\langle , \rangle$ is the Killing form on
 $\frak{g}_1$ (and $\frak{g}_2$) normalized so that
 $\langle \theta,\theta \rangle=2$ for the highest root
 $\theta$.

It is easy to see from the next  Lemma (5.2) that for a
finite dimensional representation $V$ of $\frak g_1$ given by a
Lie algebra homomorphism $\varphi : \frak{g}_1 \to sl(V)$, we have
$m_\varphi = m_V$, where  $sl(V)$
is the Lie algebra of trace $0$ endomorphisms of $V$.
\vskip1ex

 We give
 an expression for $m_V$ in the following lemma. Write the
 formal character $$ ch\, V=\sum n_{\lambda} e^{\lambda}.
 \tag"(4)"
 $$

 \vskip1ex
 \proclaim{(5.2)
 Lemma} $$ m_V=\frac{1}{2} \sum_{\lambda} n_{\lambda}
 <\lambda,{\theta}^\vee >^2.\tag"(1)"$$
  In particular, for the adjoint representation ad
 of $\frak g$ we have $$ m_{ad}=2 (1+<\rho,
{\theta}^\vee>),\tag"(2)"$$  where $\rho$
 as usual is the half sum of the positive roots of
 $\frak g$.

 Similarly, for the standard $n$-dim. representation $V^n$
  of $sl_n,
  m_{V^n}=1$.
\endproclaim
 \vskip1ex \demo{Proof}
 It suffices to show that, for the irreducible
 representation $W(m)$ (of dim $m+1)$ of $sl_2$ $$
 \frac{1}{2} \sum_{n=0}^{m} <m \rho_1-n\alpha, H>^2
=\left(\matrix m+2\\3\endmatrix\right), \tag"(3)"
$$
 where $\alpha$ is the unique positive root
 of $sl_2, H$ the corresponding coroot and $\rho_1:=
 \frac{1}{2}\alpha$. Now the left side of (3) is equal to
 $$\spreadlines{2\jot}
 \align
 2\sum_{n=0}^{m}(\frac{m}{2}-n)^2 &=4
   \sum_{k=1}^{k_o}k^2=\frac{m(m+1)(m+2)}{6},\
 \text{if}\
  m= 2 k_o\ \text{is even, and}\\
 2\sum_{n=0}^{m}(\frac{m}{2}-n)^2 &= 2\sum_{n=0}^{m}(k_o-\frac{1}{2}-n)^2 ,
 \ \text{if}\ m=2k_o-1\ \text{is odd}\\
 &=4\sum_{k=1}^{k_o} (k- \frac{1}{2})^2 =(4
   \sum_{k=1}^{k_o} k^2)+k_o-4 \sum_{k=1}^{k_o}k\\
 &=\frac{m(m+1)(m+2)}{6}\,.
 \endalign
 $$
 So in either case the
 left side of $(3)=\frac{m(m+1)(m+2)}{6}=\left(
 \smallmatrix m+2\\3\endsmallmatrix\right)$. This proves
 the first part of the lemma.

 For the assertion regarding the adjoint representation,
 we have
 $$
 ch (ad)=\text{dim} \frak h .
 e^0 + \sum_{\beta \epsilon
 \Delta_{+}} (e^{\beta}+ e^{-\beta}). $$
 So \,\,\,\,\,$m_{ad}=\sum_{\beta \epsilon \Delta_{+}}
 <\beta,{\theta}^\vee >^2$
 $$\align
 &= 4 +\sum_{\beta
 \in\Delta_+\setminus {\theta}} <\beta,\theta^\vee>
 \,\,\,\, ,\text{since} <\beta,\theta^\vee>= 0 \,\text{or} \,1, \text{for any}
\,\beta
 \in\Delta_+\setminus {\theta}\\
 &= 4+<2\rho-\theta,
 \theta^\vee>\\
 &=2(1+<\rho,{\theta}^\vee>).
 \endalign$$
 The assertion about $m_{V^n}$ is easy to
 verify.\enddemo

 \vskip1ex
\flushpar {\bf (5.3)} {\it  Remark}.  The number $(1+< \rho,{\theta}
 ^\vee>)$ is called the {\it dual Coxeter number} of
 $\frak g$. Its  value  is given as below.
 $$ \matrix \format \c\qquad &&\c\\
 \text{Type of}\ \frak g &\text{dual} &\text{
 Coxeter } &\text{number}\\ A_{\ell} &&\ell+1\\
 B_{\ell} &&2\ell-1\\ C_{\ell} &&\ell+1\\ D_{\ell}
 &&2\ell-2\\ E_6 &&12\\ E_7 &&18 \\ E_8 &&30 \\
 G_2 &&4 \\ F_4 &&9 \endmatrix
 $$

 \vskip1ex
 Now we can state the main theorem of this section.

\proclaim{(5.4) Theorem}  With the notation as in
 \S 5.1
$$ \text{Det}(\Cal U(V)) \simeq \frak L(m_V\chi_o),$$
  for any finite dimensional
 representation $V$ of $G$,
 where the number $m_V$  is  defined
 by (3) of \S 5.1.
\endproclaim
 \vskip1ex \demo{Proof}  By Proposition (C.13), there exists an integer $m$
 such that $$\text{Det}\ (\Cal U(V)) \simeq
 \frak L(m\chi_o)\in\ \text{Pic}\ ({\Cal
 G}/{\Cal P}).
  $$
 We want to prove that $m=m_V:$ Set
 $\Cal U_o:= \Cal U(V)_{\mid C\times X_o}$ as the family
 restricted to the Schubert variety $X_o :=X_{\frak s_o} $
 (cf. proof of Proposition C.13). Denote by
 $\alpha$ (resp. $\beta)$ the canonical generator of
 $H^2(X_o,\Bbb Z)$ (resp. $H^2(C,\Bbb Z )).$
  Then it suffices to
 show  that
 Det $\Cal U_o \simeq
 \frak L(m_V\chi_o)_{\mid _{X_o}},$ which is
 equivalent to showing that the first Chern class
 $$
  c_1(\text{Det}
 \,\Cal U_o)=m_V \alpha:
 \tag"(1)"
 $$

{}From the definition of the determinant bundle we have
  $$ c_1
 (\text{Det}\, \Cal U_o)= -c_1 (\pi_{2
 \ast}\Cal U_o), \tag"(2)"
$$ where $\pi_2$ is the projection $C\times
 X_o\rightarrow X_o$, and the notation $\pi_{2 \ast}$ is as
 in \cite{F, Chapter 9}.

 Since $G$ is semisimple, the associated vector bundle ${\Cal U}(V)$ has
 $$
 c_1 (\Cal U_o)=0.
 \tag"(3)"
 $$

Let
 $\tilde{\alpha}$ (resp. $\tilde{\beta})$ be  the pull
 back of $\alpha$ (resp. $\beta)$ under $\pi_2$ (resp. $\pi_1$).
 Now write
 $$ c_2
 (\Cal U_o)=l\tilde{\alpha}
 \tilde{\beta},\ \text{for some (unique)}\ l\in
 \Bbb Z. \tag"(4)"
 $$

Let $T_{\pi_2}$ be the relative tangent bundle along the fibers of
$\pi_2$. Let us denote by $c_1$ (resp. $c_2$) the first (resp. second)
Chern class of $\Cal U_o$.
By the Grothendieck's
 Riemann-Roch theorem \cite{F, \S9.1} applied to the
 (proper) map $\pi_2$, we get
 $$\align
 \text{ch}(\pi_{2\ast}
 \Cal U_o) &=\pi_{2\ast}(\text{ch}
 (\Cal U_o). \text{td}(T_{\pi_2}))\\ &= \pi_{2\ast}[(\text{rk}\, \Cal U_o
 +c_1+\frac {1}{2}(c_1^2-2c_2))(1+\frac{1}{2}c_1(T_{\pi_2}))]\\
&= \pi_{2\ast}[(\text{rk}\, \Cal U_o
 -c_2)(1+\frac{1}{2}c_1(T_{\pi_2}))],\,\,\, \text{by \,(3),}
\endalign $$
where ch denotes the Chern character and td denotes the Todd class.
Hence
 $$ \align c_1(\pi_{2\ast}
 \Cal U_o) &=\pi_{2\ast}(-c_2(\Cal U_o))\\
&= \pi_{2\ast}(-l\tilde{\alpha}
 \tilde{\beta}),\,\,\,\,\text{by \,(4)}\\
&= -l.\alpha,\,\,\,\text{since}\,\pi_{2\ast}
 (\tilde{\alpha}\tilde{\beta})=\alpha.\tag"(5)"
\endalign $$

  So to prove the theorem, by (1),(2) and (5), we
 need to show that $l=m_V$, where $l$ is given by
 (4):

 \vskip1ex It is easy to see (from its
 definition) that topologically the bundle
 $\Cal U_o$ is pull-back of the bundle
 $\Cal U'_o$ (Where $\Cal U'_o$ is the
 same as $\Cal U_o$ for $C
 ={\Bbb P}^1)$ on $ {\Bbb P}^1\times X_o$ via
 the map $$ C\times X_o \overset \delta\times \text{I}
 \to\rightarrow  {\Bbb P}^1\times X_o, $$
 where $\delta:C \rightarrow {\Bbb P}^1$ pinches
 all the points outside a small open disc around $p$ to a
 point. Of course the map $\delta$ is of degree 1, so the
 cohomology generator $\alpha$ pulls back to the
 generator $\beta $. Hence it suffices to compute the
 second Chern class of the bundle
 $\Cal U'_o$ on $
  {\Bbb P}^1\times X_o :$

 \vskip1ex Choose $X_\theta
 \in \frak g_\theta$ (where $\theta$ is
 the highest root of $\frak g)$ such that
 $<X_{\theta},-\omega X_{{\theta}} > =1$, where $\omega$ is the
 Cartan involution of $\frak g$ and $<,>$ is the
 Killing form on $\frak g$, normalized so that $<
\theta,\theta>=2$. Set $Y_\theta:=-\omega
 (X_{\theta}) \in \frak g_{-{\theta}}$.
 Define a Lie algebra homomorphism
$ sl_2 \rightarrow
   {\frak g}  \otimes_{\Bbb C} \Bbb C[t^{\pm 1}],$ by
 $$
 \matrix \format\l&&\quad\l\\
 X &\mapsto & Y_{\theta} \otimes t\\
 Y &\mapsto &
  X_{\theta} \otimes t^{-1}\\
 H &\mapsto &- \theta^\vee \otimes 1,
 \endmatrix$$
 where $\{X,Y,H\}$ is the standard
 basis of $sl_2$. The corresponding group
 homomorphism (choosing a local parameter $t$ around $p$)
 $\eta: SL_2 (\Bbb C) \rightarrow{\Cal G}$ induces a
 biregular isomorphism $\overline{\eta}:{\Bbb P}^1
\approx SL_2(\Bbb C)/B_1
 \tilde{\rightarrow} X_o$, where $B_1$
 is the standard Borel subgroup of $SL_2({\Bbb
 C})$ consisting of upper triangular matrices. In what
 follows we will identify $X_o$ with
 ${\Bbb P}^1$ under $\overline{\eta}$. The
 representation $V$ of $G$ on restriction, under the
 decomposition(2) of \S 5.1, gives rise to a continuous
 group homomorphism$$\psi : SU_2 ({\theta}) \rightarrow
 \prod_i \,(\text{Aut}  V_i),$$where $SU_2 ({\theta})$ is the
 standard compact form (induced from the involution $\omega$)
 of the group $SL_2 ({\theta})$ (with Lie algebra
 $sl_2 ({\theta}))$.

 \vskip1ex There is a principal
 $SU_2$-bundle $\Cal W$ on $S^4$ (in the
 topological category) got by the clutching construction
 from the identity map  $S^3 \approx S U_2 \rightarrow
 SU_2$. In particular, we obtain the vector bundle
 $\Cal W (\psi) \rightarrow S^4$ associated to
  the principal bundle $\Cal W$ via the representation $\psi$, which
 breaks up as a direct sum of subbundles $\Cal W_i
 (\psi)$ (got from the representations $V_i)$.

 \vskip1ex
 We further choose a degree 1 continuous map $\nu:
 {\Bbb P}^1\times {\Bbb P}^1 \rightarrow
 S^4$. We claim that the vector
 bundle $\Cal U'_o$ on
 ${\Bbb P}^1 \times {\Bbb P}^1$ is isomorphic
 (in the topological category) with the pull
 back $\nu^{\ast} (\Cal W (\psi)):$

 \vskip1ex
 Define a map $\Phi: S^1 \times (SU_2 / D)   \rightarrow
 SU_2$ by
 $$
 \left( t, \pmatrix a\,\, b \\ c\,\, d\endpmatrix\ \text{mod}\, D\right)
 \mapsto \pmatrix d\,\, ct^{-1} \\bt \, a\endpmatrix \pmatrix d\,\, c \\ b\,\,a
  \endpmatrix^{-1},
 $$
 for $\left( \smallmatrix a\, b \\ c\, d\endsmallmatrix\right)\in SU_2$ and
 $t \in S^1;$ where $D$ is the  diagonal
 subgroup of $SU_2$. It is easy to see that the
 principal $SU_2$-bundle $\nu^{\ast} (\Cal W)$
 on ${\Bbb P}^1\times{\Bbb P}^1$ is
 isomorphic with the principal $SU_2$-bundle obtained by
 the clutching construction from the map $\Phi$ (by
 covering ${\Bbb P}^1\times {\Bbb P}^1
 =S^2\times S^2=H^{+} \times S^2 \cup
 H^{-}\times S^2,$ where $H^{+}$ and $H^{-}$ are resp.
 the upper and lower closed hemispheres). By composing
 $\Phi$ with the isomorphism  $SU_2 \rightarrow SU_2
 ({\theta})$ (induced from the Lie algebra homomorphism
 $sl_2 \rightarrow sl_2({\theta})$ taking $X
 \mapsto X_{\theta},Y\mapsto Y_\theta
 $, and $H\mapsto{\theta}^\vee)$, and using the
 isomorphism $\overline{\eta}$ together with the
 definition of the vector bundle
 $\Cal U_o$ we get the assertion that
 $\Cal U'_o\approx
 \nu^{\ast}(\Cal W(\psi))$. So
 $$\align c_2
 (\Cal U'_o)&=
 \nu^{\ast} (c_2(\Cal W(\psi))) =\nu^{\ast}
 \sum_i c_2(\Cal W_i (\psi))\\
 &=\sum_i
 {\pmatrix m_i+1\\ 3\endpmatrix} \tilde{\alpha}\tilde{\beta},
 \ \text{by the following lemma}
 \ \text{(since}\ \nu \ \text{is
  a map of degree}\ 1).
 \endalign$$
 Hence $l=
 \sum_{i}\left(\smallmatrix m_i+1 \\3\endsmallmatrix\right) =m_V$, proving the
 theorem modulo the next lemma. \enddemo

 \vskip1ex
 \proclaim{(5.5) Lemma}  Let $W(m)$  be the
 $(m+1)$-dimensional  irreducible representation of $SU_2$
  and let $\Cal W(m)$  be the vector
 bundle  on $S^4$  associated to the
 principal $SU_2$-bundle $\Cal W$  on
 $S^4$  (defined in the proof of Theorem 5.4) by
 the representation $W(m)$  of $SU_2$.  Then
\endproclaim
 $$ c_2 (\Cal W(m))=\pmatrix m+2 \\ 3\endpmatrix \Omega,
 \tag"(1)"$$ {\it where $\Omega$ is the
 fundamental cohomology generator of } $H^4(S^4,
 \Bbb Z)$.
 \vskip1ex \demo{Proof} By the Clebsch
-Gordan theorem (cf.\cite{Hu, Page 126}), we have the following
 decomposition as $SU_2$-modules: $$ W(m) \otimes W (1)
=W(m+1) \oplus W(m-1) ,\ \text{ for any } m\geq 1.$$
 In particular, the Chern character $$\text{ch} \Cal W
 (m). \,\text{ch} \Cal W(1)=\text{ch} \Cal W(m+1)+\text{ch}
 \Cal W(m-1). \tag"(2)"$$ Assume,
 by induction, that (1) is true for all $l \leq
 m$. (The validity of (1) for $l=1$ is trivial to
 see.) Then by (2) we get
 $$ \matrix\format \l&&\quad\l\\
 \text{ch}\, \Cal W(m+1) &=& \text{ch} \Cal W(m).\, \text{ch}
 \Cal W(1)- \text{ch} \Cal W(m-1) \\&=&
 ((m+1).1-c_2 \Cal W(m)) (2.1-c_2
 \Cal W(1)) \\& &-(m.1-c_2 \Cal W
 (m-1)),\,\, \text{ since}\ c_1 \Cal W(l)= 0
 \ \text{ as it is a}\ SU_2\text{-bundle}.
 \endmatrix
 $$
 Hence by induction
 $$\aligned
 \text{ch} \Cal W(m+1)
 =\left( (m+1).1-
 \pmatrix m+2 \\ 3\endpmatrix \Omega\right) (2.1-\Omega) -
  &\left( m.1-
 \pmatrix m+1 \\ 3\endpmatrix \Omega\right).\\
 \endaligned
 \tag"(3)"
 $$
 Writing
 \,$\text{ch}
 \Cal W(m+1)=(m+2). 1-c_2\Cal W
 (m+1)$, and equating the coefficients from
 (3), we get
 $$\align
 c_2
 \Cal W(m+1)
 &=\left( 2\pmatrix m+2\\ 3\endpmatrix +m+1
 -\pmatrix m+1 \\ 3\endpmatrix \right) \Omega\\
 &=\pmatrix m+3 \\ 3\endpmatrix
 \Omega.
 \endalign$$
 This completes the induction and
 hence proves the lemma. \qed\enddemo
\vskip1ex

Recall that for any connected  complex simple group $G$,
the third homotopy group $\pi_3(G)$ is canonically isomorphic with
$\Bbb Z$.
\vskip1ex
\proclaim{(5.6) Corollary} For any representation $\rho$ of $G$ in
a finite dimensional vector space $V$, the induced map $\pi_3(G) \to
\pi_3(SL(V))$ is multiplication by the number $m_V$.\endproclaim

\vskip1ex \demo{Proof} The representation $\rho : G \to SL(V)$ gives
rise to a morphism $\tilde{\rho} :
{\Cal G} /{\Cal
 P} \to \Cal G^o/\Cal P^o$, where $\Cal G^o := SL(V) (\hat{\text
{{\sl k}}}_p)$ and $\Cal P^o := SL(V) (\hat {\Cal O}_p)$. Moreover, the
family $\Cal U^o(V)$ parametrized by $\Cal G^o/\Cal P^o$ (got from the
standard representation $\epsilon$ of $SL(V)$ in $V$) pulls-back to the
family $\Cal U(V)$ (parametrized by ${\Cal G} /{\Cal P}$). In particular,
from the functoriality of the determinant bundle (cf.\S 3.7), Theorem(5.4), and
Lemma (5.2), we see that the induced map  $\tilde{\rho}^{\ast} :
H^2(\Cal G^o/\Cal P^o, \Bbb Z) \to H^2(\Cal G/\Cal P, \Bbb Z)$ is
multiplication by the number $m_V$ (under the canonical identifications
$H^2(\Cal G^o/\Cal P^o, \Bbb Z) \simeq \Bbb Z \simeq  H^2(\Cal G/\Cal P, \Bbb
Z)$). But the flag variety $\Cal G/\Cal P$ is homotopic to the based loop
group $\Omega_e(K)$ (where $K$ is a compact form of $G$), and similarly
$\Cal G^o/\Cal P^o$  is homotopic to $\Omega_e(SU(V))$. In particular, by the
Hurewicz's theorem and the long exact homotopy sequence corresponding to the
fibration $\Omega_e(K) \to P(K) \to K$ (where $P(K)$ is the path space of
$K$ consisting of the paths starting at the base point $e$), the corollary
follows.\enddemo

\vskip1ex

\subheading {6. Statement of the main theorem and
 its proof}

 \vskip4ex

  \flushpar {\bf (6.1) Definitions.} Recall the definition of the homogeneous
line bundle
  $\frak L(m\chi_o)$ on $X :={\Cal G} /{\Cal
 P}$ ( for any $m \in {\Bbb Z}$)
 from  \S (C.6). Define, for any $p\in \Bbb Z$, (cf. \cite
{$\text{Ku}_1$, \S 3.8})
 $$H^p (X, \frak L (m \chi_o))= \varprojlim_{{\frak w} \in
\widetilde{W}/W}\,\,H^p(X_{\frak w},
 \frak L (m
 \chi_o)_{\mid_{X_{\frak w}}}).\tag"(1)" $$
Since any $g \in \widetilde{\Cal G}$ acts as an algebraic automorphism of the
line bundle  $\frak L(m\chi_o)$ (cf. \S C.6),
 $H^p(
X, \frak L(m
 \chi_o))$ is canonically a   $\tilde{\Cal G}$-module. This
module is determined in \cite{$\text{Ku}_1$} ( and also in \cite{M}).
We summarize the results :
$$\gather H^p (X, \frak L (m \chi_o))= 0, \,\,\text{if} \,\, p> 0\,\,
\text{and} \,\,m\geq 0, \tag"(2)"\\
H^0 (X, \frak L (m \chi_o))= 0, \,\,\text{if} \,\,m < 0,\,\text{and}
\tag"(3)"\\
H^0 (X, \frak L (m \chi_o)) \simeq L(\Bbb C,m)^\ast  \,\,\,\text{for} \, m \geq
0, \,\,\text{as}~ \widetilde{\Cal G}-\text{modules}\,,\tag"(4)"
  \endgather $$
where $  L(\Bbb C,m)$ is the integrable  highest weight (irreducible)
$\tilde{{\frak g}}$
 -module corresponding to the trivial $\frak g$-module $\Bbb C$ and the central
charge $m$ (cf. \S A.2), and $  L(\Bbb C,m)^\ast$
 denotes its full vector space dual. (By \S C.4,
 $L(\Bbb C,m)$ and hence $ L(\Bbb C,m)^\ast$ acquires a canonical structure of
$\widetilde{\Cal G}$-module.) For any subgroup $H\subset  \widetilde{\Cal G}$,
by $ H^p (X, {\frak L} (m \chi_o))^H$ we mean the space of  $H$-invariants
in $H^p (X, {\frak L} (m \chi_o))$.

 Recall the definition of the map $\varphi : \Cal G \to \Cal X_o$ from
\S 1.4, and the family $\Cal U$ parametrized by $X $ from
Proposition (2.8). Now define
$$ \align X^s &= \{ g \Cal P\in  X:
 \varphi (g)\, \text{is semistable} \} \\
&= \{x \in X: \Cal U_{\mid C\times x} \text{is semistable} \},\endalign $$
and set (for any $\frak w\in \widetilde{ W} / W)$
$$ X_{\frak w}^s =  X^s
\cap X_{\frak w}. $$
 Then by  \cite{R$_1$, Proposition (4.8)}, $ X_{\frak w}^s $ is a
Zariski open (and non-empty, since $1\in X^s_{\frak w}$) subset of $X_{\frak
w}$, in particular,
$ X^s $ is a Zariski  open  subset of $X$.
Now define
$$H^p (X^s, \frak L (m \chi_o))= \varprojlim_ { \frak w \in\widetilde{ W} /
 W}\,\,H^p(X_{\frak w}^s,
 \frak L (m
 \chi_o)_{\mid_{X_{\frak w}^s}}).\tag"(5)" $$
Clearly $\Gamma$ keeps $ X^s$ stable, and by \S 2.7,
$\Gamma$ acts morphically on the line bundle  $\frak L (m \chi_o)$ for any $m$
which is a multiple of $m_V$ (for some finite dimensional representation $V$ of
$G$), in particular,
$\Gamma$ acts on the cohomology $H^p (X^s, \frak L (m \chi_o))$, and we can
talk of the space of
 $\Gamma$-invariants $H^p (X^s, \frak L (m \chi_o))^{\Gamma}$.

The family $\Cal U_{\mid X^s}$ yields a morphism
$\psi : X^s\to \frak M$, which maps any $x\in X^s$ to the $S$-equivalence class
of the semistable bundle $\Cal U_x$, where $\frak M$ is the moduli space of
semistable $G$-bundles on $C$ (cf. Theorem 3.4).
\vskip1ex \proclaim{(6.2)
 Lemma}  There exists a $v_o
 \in \widetilde{W} / W$  such that $$ \psi
 (X_{\frak v_o}^{s})={\frak M} . $$
\endproclaim
 \demo{Proof} Since $\bigcup_{\frak w}
 X_{\frak w}^{s}= X^s$ and
 $\psi (X^s)=\frak M$, we get that
 $\frak M=\bigcup_{\frak w} \psi
 (X_{\frak w}^{s})$. But by a result of Chevalley
 (cf. \cite{B, Chapter AG, Corollary 10.2}), $\psi
 (X_{\frak w}^s)$ is a finite  union of locally closed  subvarieties
 $\{\frak M^i_{\frak w}\} $ of $\frak M$. Hence $\frak M$ is a
countable union $\bigcup \frak M^i_{\frak w}$ of  locally closed
subvarieties. But then, by a Baire category argument, $\frak M$ is a certain
finite union
of (locally closed) subvarieties $\{\frak M^1_{\frak w_1},\dots , \frak
M^n_{\frak w_n} \}$. Now choosing a
  $\frak v_o  \in \widetilde{W} / W$  such that $\frak v_o \geq  \frak w_i$ ,
 for all  $1\leq i \leq n$, we get that
 $ {\frak M} = \psi
 (X_{\frak v_o}^{s})$. This proves the
 lemma.\enddemo
\vskip1ex \proclaim{(6.3) Corollary} The moduli space
$\frak M$ is a unirational variety.
\endproclaim
\demo{Proof} Since $X_{\frak w}^{s}$ is an open subset of
 $X_{\frak w}$ and $X_{\frak w}$ is a rational variety
( by the Bruhat decomposition), the corollary
follows from the above lemmma (6.2). \enddemo
 \vskip1ex

\proclaim{(6.4) Proposition}  For any $d \geq 0$ and any finite dimensional
representation $V$ of $G$, the canonical
 map $$\psi^{\ast}:H^0({\frak M},\, \Theta
  (V)^{\otimes d})\rightarrow H^0
 (X^s,\psi^{\ast}(\Theta
  (V))^{\otimes d})^{\Gamma}$$
is an isomorphism, where $\Theta (V)$ is the theta bundle on the
moduli space $\frak M$ associated to the representation $V$ (cf. \S3.7),
and the vector space on the right denotes the space of $\Gamma$-invariants
under its natural action on the line bundle $\psi^{\ast}(\Theta (V))$. (Since
the map
$\psi: X^s \to \frak M$ is $\Gamma$-equivariant, with trivial action
of $\Gamma$ on $\frak M$, the pull-back bundle $\psi^{\ast}(\Theta (V))$
has a natural $\Gamma$-action.)
 \endproclaim
\demo{Proof} Using Lemma (6.2) we see that the map $\psi^{\ast}$ is injective.
Now the second
part of Proposition (2.8), and Proposition (3.6) show that if $x$ and $y$
are two points in $X^s$ with $\Cal U_y \simeq \text{gr}(\Cal U_x)$, then
$y$ belongs to the Zariski closure of the $\Gamma$-orbit of $x$. In
particular, two points in $X^s$ are  in the same fiber of $\psi$
if and only if the closures of their $\Gamma$-orbits intersect. This,
in turn, shows that if $\sigma$ is a $\Gamma$-invariant regular section of
$\psi^{\ast}(\Theta
  (V))^{\otimes d}$ on $X^s$, it is induced from a set theoretic section
$\underline
{\sigma}$ of $\Theta (V)^{\otimes d}$ on $\frak M$. That $\underline
{\sigma}$ is regular, is seen by taking any   Schubert variety
 $X_{\frak w}$ such that $\psi(X_{\frak w}^s)= \frak M$ (cf. Lemma 6.2) and
applying
Proposition (4.1) to the morphism $\psi_{\mid X^s_{\frak w}} : X^s_{\frak w}\to
\frak M$. \qed \enddemo
 \vskip1ex
 By the functorial property of the theta bundle, $\Theta (\Cal U(V))_{\mid
X^s}$ is canonically
isomorphic to $\psi^{\ast}(\Theta (V))$, since $\psi$ is defined using the
restriction of the family $\Cal U(V)$ to $X^s$ (cf. \S3.7). Moreover, as
observed in \S3.7, the line bundles $\Theta (\Cal U(V))$ and $\text{Det} (\Cal
U(V))$ coincide on the whole of $X$.

\proclaim{(6.5) Proposition} Any $\Gamma$-invariant regular section
of $\psi^{\ast}(\Theta
  (V))^{\otimes d}$ on $X^s$ extends uniquely to a regular section of (Det
$\Cal U(V))
^{\otimes d}$ on $X$.
\endproclaim
This proposition will be proved in the next section.
\vskip1ex
 We now state and prove our main theorem, assuming the validity of Proposition
(6.5).
\vskip1ex
\proclaim{(6.6) Theorem} Let
 the triple $\frak T= (G,C,p)$ be as in \S1.1, and let $V$ be a
 finite dimensional  representation
 of $G$.  Then,  for any $d \geq 0$,
 $$
 H^0({\frak M},\,\Theta
  (V)^{\otimes d}) \simeq H^0 ({\Cal G} /
{\Cal P},\frak L (dm_V\chi_o))^{\Gamma},$$
where the latter space of $\Gamma$-invariants is defined in
\S 6.1,  the integer $m_V$ is the Dynkin index of $V$ defined in \S5.1, and
the moduli space $\frak M$ and the theta bundle $\Theta (V)$ are
as in \S\S3.4 and 3.7 respectively.

In particular, $H^0 ({\Cal
 G} /{\Cal P}, \frak L (dm_V\chi_o))^{\Gamma} $
 is finite dimensional.

 (Observe that by (4) of
\S 6.1,  $H^0 ({\Cal
 G} /{\Cal P}, \frak L (dm_V\chi_o))^{\Gamma} $ is isomorphic with
the space of $\Gamma$-invariants in the dual space
$L(\Bbb C,dm_V)^\ast$.)
\endproclaim
\demo{Proof} We first begin with some simple observations:
\vskip1ex
(a) {\it For any algebraic line bundle $\frak L $ on $X$, the canonical
restriction map
$H^0(X, \frak L) \to H^0(X^s, \frak L_{\mid X^s})$ is injective}:
This is seen by restricting a section to each Schubert variety
$X_{\frak w}$, and observing that
$X_{\frak w}^s$ is non-empty,
 and open (and hence
dense) in the irreducible variety $X_{\frak w}$.

(b) {\it If $\frak L$ is a $\Gamma$-equivariant line bundle
on $X$ (with respect to the standard action of $\Gamma$ on $X$) (cf. \S B.7)
and
$\sigma$ is a regular section of $\frak L$ such that its restriction to
$X^s$ is $\Gamma$-invariant, then $\sigma$ itself is $\Gamma$-invariant}:
By $\Gamma$-invariance, for $\gamma \in \Gamma$, the section $\gamma^{\ast}
(\sigma)-\sigma$
vanishes on $X^s$ (and hence on the whole of $X$).

(c) {\it Suppose that $\frak L'$ and $\frak L''$ are two $\Gamma$-equivariant
line bundles on $X^s$. Then any biregular isomorphism of line bundles
$\xi : \frak L' \to
\frak L''$ (inducing the identity on the base) in fact is $\Gamma$-equivariant.
In particular, $\xi$ induces an isomorphism of the corresponding spaces of
$\Gamma$-invariant regular sections}:

Define a map $\epsilon:
 \Gamma \times X^s\rightarrow
 \Bbb C^{\ast}$ by
$$ \epsilon (\gamma,x)= L_{\gamma^{-1}}\,\xi_{\gamma x}
 \,L_{\gamma}\, (\xi_x)^{-1}
 \in \, \text{Aut}_{\Bbb C} (\frak L''_x)
 =\Bbb C^{\ast}, $$
for $\gamma
 \in \Gamma$ and $x \in X^s$,
where $L_\gamma$ is the action of $\gamma$ on the appropriate line bundles,
and $\xi_ x$ denotes the restriction of $\xi$ to the fiber over
$x \in X^s$.
  It is easy to see that
 $\epsilon$ is a regular map, and  of course
 $\epsilon (1, x) = 1$  for all  $x \in X^s$.
In particular, by Corollary (2.4), $ \epsilon (\gamma, x)= 1$,
for all $\gamma \in \Gamma$. This proves  assertion (c).
\vskip1.5ex

We now consider (Det $ \Cal U(V))^{\otimes d}_{\mid X^s}$ as a
$\Gamma$-equivariant line bundle by transporting the natural $\Gamma$-action
on $\psi^{\ast}(\Theta(V))^{\otimes d}$ (cf. Proposition 6.4), via the
canonical identification
$$
 \text{Det}\, \Cal U(V)_{\mid X^s} \simeq \psi^{\ast}(\Theta(V)).
\tag"(1)"$$

Choose an isomorphism of line bundles on $X$
$$ \xi_o : (\text{Det}\, \Cal U(V))^{\otimes d} \to \frak L(\chi_o)^{\otimes
dm_V}\,,$$
 which exists by Theorem (5.4). Recall from \S 2.7  that
 $\frak L(\chi_o)^{\otimes dm_V}$ is a  $\Gamma$-equivariant
line bundle on $X$. Hence by (c) above, the map $\xi :={\xi_o}_{\mid X^s}$ is
automatically $\Gamma$-equivariant.
We have the following commutative diagram:
 $$
\CD
 H^0(X, (\text{Det}\, \Cal U(V))^{\otimes d}) @>\underset\sim\to
 {\overline{\xi_o}}>> H^0(X, \frak L(\chi_o)^{\otimes dm_V})  \\
 @VVV        @VVV\\
 H^0(X^s, (\text{Det}\, \Cal U(V))^{\otimes d}) @>\underset\sim\to
{\overline{\xi}}>>  H^0(X^s, \frak L(\chi_o)^{\otimes dm_V})\,
 \endCD
$$
where $\overline{\xi}$ (resp. $\overline{\xi_o}$) is induced from $\xi$ (resp.
$\xi_o$), and the vertical maps are the canonical restriction maps. Observe
that $\overline{\xi}$ is $\Gamma$-equivariant (since $\xi$ is so).

Further we have
$$ \align
H^0(\frak M, \Theta(V)^{\otimes d}) &\simeq  H^0(X^s, (\text{Det}
\,\Cal U(V))^{\otimes d})^{\Gamma}    \,\,\,\,\,(\text{by (1) and Proposition
6.4})\\
& \simeq H^0(X^s, \frak L(\chi_o)^{\otimes dm_V})^{\Gamma}
\,\,\,\,(\text{under} \,
\overline{\xi}).
\endalign $$

We complete the proof of the theorem by showing that the restriction map
$$ H^0(X, \frak L(\chi_o)^{\otimes dm_V})^{\Gamma} \to
H^0(X^s, \frak L(\chi_o)^{\otimes dm_V})^{\Gamma} $$
is an isomorphism:

It suffices to show that any $\Gamma$- invariant
section $\sigma$ of $\frak L(\chi_o)^{\otimes dm_V}$ over $X^s$
extends to a section over $X$, for then the extension will automatically
be $\Gamma$-invariant by (b) and unique by (a). By the above commutative
diagram,
this is equivalent to showing that any $\Gamma$-invariant section
$\sigma_o$ of (Det$ \,\Cal U(V))^{\otimes d}$ over $X^s$ extends to the whole
of
$X$. But this is the content of Proposition (6.5), thereby completing the
proof of the theorem. \qed \enddemo

\proclaim{ (6.7) Proposition}  For any $d \geq 0$, and finite dimensional
representation $V$ of $G$,  we have

$$ [L(\Bbb C,dm_V)^\ast]^\Gamma = [L(\Bbb C,dm_V)^\ast]^{\text{Lie}\, \Gamma} =
[L(\Bbb C,dm_V)^\ast]^{\frak g \otimes k[C^\ast]},$$
where $L(\Bbb C,dm_V)$ is canonically an algebraic $\Gamma$-module as in \S2.7,
$\frak g$ is the Lie algebra of the group $G$ and (as in \S1.1)
 $k[C^\ast]$ is the ring of regular functions on the affine curve $C^\ast$.

\endproclaim

\demo{Proof} Abbreviate $L(\Bbb C,dm_V)$ by $V$. Fix $v \in V$ and consider the
morphism $\pi_v : \Gamma \to V$ given by $\pi_v(\gamma) = \gamma.v$ for
$\gamma \in \Gamma$. Recall (cf. Lemma B.6) that, by definition, the action of
the Lie
algebra Lie $\Gamma$ on $ v \in V$ is given by the induced map
$(d\pi_v)_e : T_e(\Gamma)=$ Lie $\Gamma \to T_v(V) = V$.

Fix  $\theta \in V^\ast$. For any $v \in V$, define the map
$\theta_v : \Gamma \to \Bbb A^1$ by $\theta_v(\gamma) = \theta (\gamma.v)$. The
induced map $(d\theta_v)_e : T_e(\Gamma)= $ Lie $\Gamma \to T_{\theta (v)}(\Bbb
A^1) = \Bbb A^1$ is given by
$$  (d\theta_v)_e (a) = \theta (a.v), \,\,\text{for}\,
a \in \, \text{Lie} \,\Gamma. \tag"(1)"$$

 For any $\gamma_o \in \Gamma$, we now determine the map
$(d\theta_v)_{\gamma_o}$:
Consider the right translation map $R_{\gamma_o} : \Gamma \to \Gamma,$ given
by  $R_{\gamma_o}(\gamma) = \gamma \gamma_o.$ Then we have
$$  (d\theta_v)_{\gamma_o} o (dR_{\gamma_o})_e =
 (d\theta_{\gamma_o.v})_e. \tag"(2)"$$

If $\theta \in [V^\ast]^\Gamma,$ then  $\theta_v$ (for any fixed $v \in V$)
is the constant map $\gamma \mapsto \theta(v).$ In particular, $(d\theta_v)_e
\equiv 0$, proving (by 1) that $\theta \in [V^\ast]^{\text{Lie}\, \Gamma}$.
Conversely, take  $\theta \in [V^\ast]^{\text{Lie}\, \Gamma}$. Then by (1)
and (2), for any fixed $v \in V , (d\theta_v)_{\gamma_o} \equiv 0$
for any $\gamma_o \in \Gamma$. In particular,
for any fixed $v \in V$ and  $i \geq 0$,  the map $\theta_{v_ {\mid \Gamma_i}}
:
\Gamma_i \to \Bbb A^1$ ($\Gamma_i$ is as in proof of Lemma 2.1) is constant on
the irreducible components of $\Gamma_i$ (as the base field is of char. 0). But
since $\Gamma $ is connected (cf. proof of Corollary 2.4), $\theta_v$ itself is
forced to be a
constant. Thus, we have $(\gamma \theta - \theta) v = 0$, for every $v \in V$
and $\gamma \in \Gamma$;  proving  that  $\theta \in [V^\ast]^\Gamma.$ Finally,
by Lemma (2.2), we have Lie $\Gamma = \frak g \otimes k[C^\ast]$. This proves
the proposition. \qed \enddemo
\vskip5ex
\flushpar {\bf (6.8)} {\it Remarks.} (a) From the proof it is clear that the
above proposition is true with $L(\Bbb C,dm_V)$ replaced by any algebraic
representation
of the algebraic group $\Gamma$.

(b) In Conformal Field Theory, the space of vacua is defined to be the space
of invariants $[L(\Bbb C,d)^\ast]^{\frak g \otimes k[C^\ast]}$ of the
Lie algebra $\frak g \otimes k[C^\ast]$ (cf. \cite{TUY, Definition 2.2.2}).
We see, by Theorem (6.6) and Proposition (6.7), that
the space of vacua is isomorphic to the space of generalised theta functions.

 (c)  Assertion (c) in the proof of
 Theorem (6.6) is true with $X^s$ replaced by $X$.  (In
 fact, in this case we do not even need to use Corollary
 (2.4), but need only the connectedness of $\Gamma$.)  We
 outline an argument:

 Following the same notation as in the proof of assertion
 (c), in this case, for any fixed $\gamma \in \Gamma $
 the map $\epsilon_{\vert_{\gamma \times X}}:\gamma
 \times X\to \Bbb C^\ast$ is a constant $\alpha_{\gamma }$
 (since $X$ is a connected projective ind-variety).  From
 this it can been easily seen that the map $\alpha:\Gamma
 \to \Bbb C^\ast$ taking $\gamma \mapsto \alpha_\gamma $
 is a group morphism.  In particular, the derivative
 $d\alpha:$ Lie $\Gamma\to\Bbb C$ is a Lie algebra
 homomorphism.  Since the commutator [Lie $\Gamma$, Lie
 $\Gamma ] =$ Lie $\Gamma $, we get that $d\alpha\equiv
 0$.  Hence, by an argument used in the proof of
 Proposition (6.7), we see that the map
 $\alpha$ itself is identically 1.  This proves
 assertion (c) for $X$.

 \vskip1ex

 Now, as shown in \cite{KN}, the algebraic
 $\Gamma$-action on the line bundle (Det $\Cal U(V)$)$^{\otimes
 d}_{\vert_{X_s}}$ can be extended to an algebraic
 $\Gamma$-action on the line bundle (Det  $\Cal U(V)$)$^{\otimes
 d}\to X$.  In particular, the isomorphism $\xi_0:$ (Det
 $\Cal U(V)$)$^{\otimes d} \to \Cal L(\chi_0)^{\otimes dm_V}$ (cf.
 proof of Theorem 6.6) is $\Gamma$-equivariant and hence
 so is its restriction to $X_s$.  This provides a proof
 of Theorem (6.6), without using Corollary 2.4 (but using
 only the connectedness of $\Gamma$).

 \vskip1ex

As an immediate consequence of the above remark (b), we obtain the
following.

\vskip1ex
\proclaim{(6.9) Corollary}  Let
 the notation and assumptions be as in
 Theorem (6.6).  Then  the space
 of covariants $L(\Bbb C,dm_V) / ((\frak g
 \otimes_k  k[C^\ast] ). L(\Bbb C,dm_V))$
   is finite dimensional. (Cf. \cite{K, Exercise 11.10, p. 209} for a purely
algebraic proof of this corollary.)
 \endproclaim
\vskip6ex

 \subheading {7. Proof of Proposition (6.5)}
\vskip4ex

\proclaim{(7.1) Lemma}  Let $X$ be an irreducible normal  variety,\,\, $U
\subset X$ a non-empty open subset and $\frak L$ a line bundle on $X$.  Then
any
element of $\underset {n\in \Bbb Z_+}\to\oplus H^0(U,\frak L^n)$ which is
integral over $\underset
n\to\oplus H^0(X,\frak L^n)$ belongs to $\underset n\to\oplus H^0(X,\frak
L^n)$.
\endproclaim
\demo{Proof} Since the  rings in question are graded, it suffices to
prove the lemma only for homogeneous elements. Let $b\in H^0(U,\frak L^{n_o})$
be integral over $\oplus H^0(X,\frak L^n)$, i.e., $b$ satisfies a
 relation $b^m+a_1 b^{m-1}+\dots + a_m=0$ with $a_i\in \oplus
H^0(X,\frak L^n)$.  Let $D$ be a prime divisor in $X\setminus U$ and let $b$
have a pole of order
$\ell \geq 0$ along  $D$.  Then the order of the pole of $b^m$ along $D$ is of
course $\ell m$ and that
of $a_ib^{m-i}$ is $\leq \ell (m-1)$ for every $i\geq 1$.  But since
$b^m+a_1b^{m-1}+\dots + a_{m-1}b$ \,\, is by assumption regular along $D$, we
 are forced to have  $\ell =0$, i.e., $b$ is regular
along $D$.  Hence $b\in H^0(X,\frak L^{n_o})$. \qed \enddemo

 \proclaim{(7.2) Lemma}  Let $f:X\to Y$ be a morphism
 between projective varieties and $\Cal L$ an ample line
 bundle on $Y$.  Then the ring $\underset n\geq
 0\to\oplus H^{0}(X,f^\ast \Cal L^n)$ is integral over
 the ring $\underset{n\geq 0}\to{\oplus}H^{0}(Y,\Cal L^n)$.
 \endproclaim

 \demo{Proof} First of all,
$$\underset n\geq 0\to\oplus
 H^{0}(X,f^\ast \Cal L^n) \approx H^{0} (Y,\Cal L^n
 \otimes_{\Cal O} f_\ast \Cal O)\,. \tag"(1)" $$
  Since $f_\ast\Cal O $ is
 coherent, we can write
 $$
 0 \to \Cal K \to \Cal L_1\oplus \dots \oplus \Cal L_m
 \to
 f_\ast\Cal O  \to 0,
 $$
 for some line bundles $\Cal L_i$ on $Y$.
 Since $\Cal L$ is ample on $Y$, there exists $n_o \geq 1$ such
 that for $n\geq n_o$, $H^1 (Y,\Cal L^n\otimes \Cal K)=0$;
 in particular,
 $$H^{0} (Y,\Cal L^n \underset{\Cal O
 }\to{\otimes} (\underset {i}\to{\oplus}\Cal L_i) )
 \twoheadrightarrow H^0 (Y,\Cal L^n \underset{\Cal O
 }\to{\otimes} f_\ast \Cal O) \tag"(2)"$$
 is surjective for $n\geq
 n_o$.
  We now prove that $\underset {n\geq 0}\to{\oplus} H^0
 (Y,\Cal L^n \otimes \Cal L_i)$ is finitely generated
 over $\underset {n\geq 0}\to{\oplus} H^0 (Y,\Cal L^n), $
 for any line bundle $\Cal L_i$ on $Y$ :

 Consider the sheaf exact sequence (where $\Delta (Y)
 \subset Y\times Y$ is the diagonal)
 $$
 0 \to \Cal I_{\Delta (Y)}
 \to \Cal O_{Y\times Y} \to
 \Cal O_Y \to 0.
 $$
 There exist $\ell _0, m_0\geq 1$ such that for $n\geq
 \ell _0$ and $m\geq m_0$, $H^1(Y\times Y, \Cal I _{\Delta (Y)}
 \otimes  (\Cal L^n \boxtimes (\Cal L^m\otimes \Cal
 L_i))) =0$.  In particular, for $n\geq \ell _0$ ,
 $$
 H^{0} (Y,\Cal L^n) \otimes
 H^0 (Y,\Cal L^{m_0} \otimes \Cal L_i)
 \twoheadrightarrow
 H^0 (Y,\Cal L^{n+m_0} \otimes \Cal L_i)
 $$
 is surjective.  This, in particular (using (2)), proves that the
 ring $\underset{n\geq 0}\to{\oplus} H^0 (Y,\Cal L^n
 \underset{\Cal O}\to{\otimes} f_\ast(\Cal O))$ is
 finitely generated over $\underset{n\geq 0}\to{\oplus}
 H^{0} (Y,\Cal L^n)$, in particular, is integral (cf. \cite{AM}).  This proves
the lemma by (1).\qed
 \enddemo

 \proclaim{(7.3) Proposition} Let $Y$ be a normal
 irreducible variety parametrizing a family $\Cal E $ of
 $G$-bundles.  Consider the induced map $\beta : Y^s\to
 {\frak M}$ (cf. Theorem 3.4; where $Y^s$ is the subset of $Y$ consisting of
those
 $y\in Y$ such that $\Cal E _y$ is semistable).  Fix
 a representation $V$ of $G$ and fix an isomorphism
 $$
 \Theta (\Cal E (V))_{\vert_{Y^s}}
 \approx \beta ^\ast (\Theta(V)).
 $$

  Then for any section $\sigma \in
 H^0(\frak M,\Theta (V)^{\otimes d})$ ($d\geq 0$), the pull-back
 section $\beta ^\ast \sigma $ extends to an element of
 $H^0(Y, \Theta(\Cal E (V))^{\otimes d})$.
 \endproclaim

 \demo{Proof} Consider the diagram (cf. \S3.7)
$$\matrix
\Theta (\Cal E (V))^{\otimes d}_{\vert_{Y^s}}
 &\longerrightarrow &\Theta (V) ^{\otimes d} &\longerrightarrow
  &\Theta ^{\otimes d}\\  \vspace{1\jot}
   @VVV @VVV  @VVV\\ \vspace{2\jot}
 Y^s &\underset\beta\to{\longerrightarrow} &\frak M &\longerrightarrow
   &\frak M_0.
\endmatrix
$$

 By Lemma (7.1), it suffices to show that $\beta ^\ast
 \sigma \in H^0(Y^s,\Theta(\Cal E (V))^{\otimes d})$ is
 integral over $\underset n\geq 0\to\oplus H^{0} (Y, \Theta
 (\Cal E (V))^{\otimes n})$:

 Since $\sigma $ is integral over $\underset n\geq
 0\to\Theta H^0(\frak M_o,\Theta^{\otimes n})$
 (by  Lemma 7.2), there exists a relation
 $$
 \sigma ^p + a_1\sigma ^{p-1} +\dots + a_p=0, \qquad
 \text{for some } a_i \in \underset {n\geq 0}\to{\oplus}
 H^0(\frak M_o, \Theta^{\otimes n}).
 \tag"($\ast$)"
 $$
 Now we have
 $$
  \matrix
 \format \c&&\quad\c\\
  Y^s &\hookrightarrow  &Y ^s_o \\\vspace{2\jot}
  \sideset \beta\and\to\downarrow
 &&\sideset\and\beta_o\to\downarrow \\\vspace{2\jot}
  \frak M &\rightarrow & \frak M_o
 \endmatrix
 $$
 where $Y^s_o := \{ y\in Y: \Cal E _y(V) \text{ is a
 semistable $GL(V)$-bundle}\}$.    Assume that $\beta
 _o^\ast a_i$ can be extended to the whole of $Y$ (which is possible by  the
 following Lemma 7.5), then we get from ($\ast $)
 $$
 (\beta ^\ast \sigma)^p + (\beta ^\ast _0
 a_1)_{\vert_{Y^s}}
 (\beta ^\ast \sigma)^{p-1} +\dots +
 (\beta ^\ast _0 a_p) _{\vert_{Y^s}} =0,
 $$
 which proves the proposition. \qed
 \enddemo

\noindent {\bf (7.4)  G.I.T. construction of the  moduli of
 vector bundles.}
 We recall the construction of the moduli space of vector
 bundles on $C$ using G.I.T.. Let $r\geq
 1$ and $\delta$ be integers. For the  fixed  point $p\in C$ and
 for a coherent sheaf {\bf F} on $C$, put {\bf F}$(m)= ${\bf F}$\otimes_{\Cal
O}
 \Cal O(mp)$, for any $m\in \Bbb Z$, where $\Cal O= \Cal O_C$ is the structure
sheaf of $C$. We can choose  a large enough
 integer $ m_o=m_o(r,\delta)$ such that for
 any semistable vector bundle {\bf E} of rank $r$ and degree
 $\delta$ on $C$, we have $H^1(\text{\bf E}(m_o))=0$ and {\bf E}$(m_o )$ is
 generated by its global sections. Let $q=\dim H^0(\text{\bf E}(m_o)) =
 \delta+r(m_o+1-g)$ (where $g$ is the genus of $C$) and
 consider the {\it Grothendieck quot  scheme} $Q$ consisting of
 coherent sheaves on $C$ which are quotients of $\Bbb C^q\otimes_{\Bbb C} \Cal
O$ with
 Hilbert polynomial (in the indeterminate $v$) $rv+q$ .
 The group $GL(q,\Bbb C)$ operates canonically on $Q$ and the action
 on $C\times Q$ (with the trivial action on $C$) lifts  to an
 action of the tautological sheaf $\Cal T $ on $C\times
 Q$.

 We denote by $R_o$ the $GL(q)$-invariant open subset of
 $Q$ consisting of those $x\in Q$ such that $\Cal T _x=\Cal T_{\mid
 C\times x}$ is locally free and such that the following
canonical map is an isomorphism:
 $$
 \Bbb C^q=H^0(\Bbb C^q\otimes_{\Bbb C} \Cal O) \tilde{\rightarrow}
 H^0(\Cal T_x).
 $$

   Then $R_o$ is smooth and irreducible and $\Cal T_{\vert C \times R_o}$ is a
rank-$r$ vector bundle.  Define the open subset (of $R_o$) $R_o^s = \{x\in R_o
: \Cal T_x \text{~is ~semistable~} \}.$
If we choose sufficiently large $m_o$ ,
 the
 G.I.T. quotient $R_o^s// GL(q)$ yields
 the {\it moduli space $\frak M_o$ of vector bundles} of rank $r$ and
 degree $\delta$. (We choose such a $m_o$ in the sequel.)
 (For all this, see \cite{NRa, Appendix A} or \cite{Le}.)

Now  let $\Cal V_o\rightarrow C\times T_o$ be a family of
 vector bundles of rank $r$ and degree $\delta$
 (parametrized by a variety $T_o$). We can find an integer
 $m_{T_o}$ such that for $m\geq m_{T_o}$, we have :
 \roster
 \item $R^1p_{2_\ast}(\Cal V_o(m))=0$.
 \item $p_{2_\ast} (\Cal V_o(m))$ is a vector bundle on $T_o$ ( of rank
 $q:= \delta+r(m+1-g))$, and
 \item the canonical map
 $p_{2}^\ast p_{2_\ast}(\Cal V_o(m))\rightarrow \Cal V_o(m)$ is surjective,
 \endroster

where $p_2: C\times T_o \to T_o$ (resp. $p_1 : C\times T_o \to C$) is the
projection on the second  (resp first) factor,  and $\Cal V_o(m)$ := $\Cal
V_o\otimes_{\Cal O_{C\times T_o}}p^\ast_1 \Cal O(mp)$.

 \vskip1ex

 Choose $\overline{m}_o$ larger than $m_{T_o}$ and $m_o$, where $m_o$ is as
above. Let {\bf P}$_o$ be the frame
 bundle of $p_{2_\ast} (\Cal V_o(\overline{m}_o))$ with the projection  $\pi_o$
: \,{\bf P}$_o\rightarrow T_o$.  Then there exists a canonical
$GL(q)$-equivariant
  morphism $\varphi _o$: {\bf P}$_o\rightarrow R_o$
 such that the families $\pi^\ast_o (\Cal V_o)$ and $\varphi
 _o^\ast (\Cal T (-\overline{m}_o))$ are isomorphic, where $\Cal T
 (-m):= \Cal T \otimes_{\Cal O_{C\times R_o}} {p}^\ast _C\Cal O (-mp)$ and
${p}_C : C\times R_o \to C$
is the projection on the first factor.

\vskip1ex

With the same notation as in Proposition (7.3) and its proof,
we have the following:

 \proclaim{(7.5) Lemma}  For any $\sigma _o\in H^0
 (\frak M_o,\Theta^{\otimes d})$ $(d\geq 0)$ the
 pull-back section $\beta^\ast _o\sigma _o\in
 H^0(Y^s_o,\Theta ({\Cal E}(V))^{\otimes d})$ can be
 extended to the whole of $Y$.
 \endproclaim

 \demo{Proof} In the  construction as in \S7.4, take $r=\dim V,\
 \delta=0$, and $T_o=Y$.

 Consider the diagram (where the map $\pi$ is the quotient map and the other
maps are as explained in
 \S\S7.3-7.4):
 $$
   \gather
     P_o\\
  \sideset\pi_o\and\to\swarrow
 \qquad\qquad
  \sideset\and\varphi_o\to\searrow \\
  \vspace{2\jot}
   Y^s_o \hookrightarrow Y \qquad\qquad R_o
 \hookleftarrow R^s_o\\
 \vspace{2\jot}
   \sideset\beta_o\and\to\searrow
  \qquad\qquad\qquad
  \sideset\and\pi\to\swarrow\\
  \vspace{2\jot}
   {\frak M}_o
   \endgather
  $$

 Now since $\text{codim}_{R_o}(R_o\backslash R_o^s)\geq 2$  and $R_o$
 is smooth (in particular normal), the section $\pi^\ast \sigma _o$
 extends (say to $\overline{\pi^\ast \sigma _o}$) to
 the whole of $R_o$.  Pull $\overline{\pi^\ast \sigma
 _o}$ via the GL$(q)$-equivariant morphism $\varphi _o$, and then push  via
the GL$(q)$-bundle  $\pi
 _o$.  This  gives an extension of the section
 $\beta_o^\ast \sigma _o$ to the whole of $Y$.  This
 proves the lemma, thereby proving Proposition (7.3)
 completely. \qed
 \enddemo
\vskip2ex
 Finally we prove Proposition (6.5) and thus complete the
 proof of Theorem (6.6).
 \vskip2ex
\noindent {\bf (7.6)} {\it Proof of Proposition (6.5).}
 Let $\tilde{\sigma}$ be a $\Gamma$-invariant section of
 $\psi^\ast (\Theta(V))^{\otimes d}$ on $X^s$. By Proposition (6.4),
 there is a section $\sigma$ of $\Theta(V)^{\otimes d}$
 over $\frak M$ such that $\psi^\ast (\sigma)=\tilde{\sigma}$.
 Let $X_{\frak w}$ be a Schubert variety.  Since $X_{\frak w}$ is irreducible
and normal (cf. \cite{Ku$_1$, Theorem 2.16}), by Proposition (7.3),
$\tilde{\sigma}
 _{\mid X^s_{\frak w}}$ extends to a (unique) section $\overline{\sigma}_\frak
w$ of  $(\text{Det } \Cal U(V))^{\otimes d}$ on $X_\frak w$.  By the uniqueness
of extensions, it is clear that for any $\frak v \leq \frak w ,
\overline{\sigma}_ {\frak w_{\vert X_\frak v} } = \overline{\sigma}_\frak v $.
In particular,
the sections $\{\overline{\sigma}_\frak w \}$ give rise to a section
$\overline{\sigma}$ of   $(\text{Det } \Cal U(V))^{\otimes d}$ on the whole of
$X$, extending the section $\tilde{\sigma}$ . This completes the proof of the
proposition.
\qed
\vskip4ex

\flushpar
  {\bf Appendix A. Affine Lie algebras - Basic definitions and their
 representations.}
\vskip 5mm
 The basic reference for this section is Kac's book [K].
 \vskip 3mm
 \flushpar
  {\bf (A.1) Definition}. Let $\frak{g}$ be a finite
 dimensional complex simple Lie algebra.  (We also fix a  Borel
 subalgebra $\frak{b}$ and a Cartan subalgebra
 $\frak{h}\subset \frak{b}$ of ${\frak g}$.) Then the associated {\it affine
Kac-Moody Lie
 algebra} is by definition  the
  space $$\tilde{\frak{g}}:= \frak{g}
 \otimes_{{\Bbb C}} \Bbb C [t^{\pm 1}] \oplus \Bbb C K$$ together with  the Lie
bracket (for $X, Y \in \frak{g}$ and $P,Q \in \Bbb C [t^{\pm 1}]$)
 $$
 [X \otimes P, Y \otimes Q] = [X,Y] \otimes PQ+(<X,Y>
 \underset t=0 \to{\text{Res}} (\frac{dP}{dt} Q))K, \, \text{~and~}
 $$
 $$[\tilde{\frak{g}}, K]=0~, $$
where $<,>$ is the
 Killing term on $\frak{g},$ normalized so that $< \theta,
 \theta> =2$ for the highest root $\theta$ of $\frak{g}$. We also define  a
certain completion
$\tilde\frak{g}_{\text{comp}}$ of $\tilde\frak{g}$
by
$$\tilde\frak{g}_{\text{comp}} = \frak{g}
 \otimes_{{\Bbb C}} \Bbb C ((t)) \oplus \Bbb C K ~,$$
where ${\Bbb C}((t))$ is the field of Laurent power series. Then
$\tilde\frak{g}_{\text{comp}}$ is a Lie algebra under the same bracket as
above.

 The Lie algebra  $\frak{g}$ sits as a Lie subalgebra of
 $\tilde{\frak{g}}$ as $\frak{g} \otimes t^0$. The Lie
 algebra $\tilde{\frak{g}}$ admits a distinguished
 `parabolic' subalgebra
 $$\tilde{\frak{p}}:= \frak{g} \otimes \Bbb C [t] \oplus
 \Bbb C K.$$
 We also define its `nil-radical' $\tilde{\frak{u}}$ (which
 is an ideal of $\tilde{\frak{p}}$) by
 $$\tilde{\frak{u}}:= \frak{g} \otimes t \Bbb C [t],$$
 and its `Levi component' (which is a Lie subalgebra of
 $\tilde{\frak{p}})$
 $$\tilde{\frak{p}}^o:= \frak{g} \otimes t^o \oplus \Bbb C
 K.$$
 Clearly (as a vector space)
 $$\tilde{\frak{p}}= \tilde{\frak{u}} \oplus
 \tilde{\frak{p}}^o.$$

Define the loop algebra $L({\frak g}) := {\frak g}\otimes_{\Bbb C} {\Bbb
C}[t^{\pm 1}]$ with Lie bracket $
 [X \otimes P, Y \otimes Q] = [X,Y] \otimes PQ $, for  $X, Y \in {\frak g}$ and
$P,Q \in \Bbb C [t^{\pm 1}].$ Then $\tilde{\frak g}$ can be viewed as a
one-dimensional central extension of
$L({\frak g})$:
$$ 0 \to {\Bbb C}K \to \tilde{\frak g} \overset{\pi}\to \longrightarrow
L({\frak{g}}) \to 0,$$
where the Lie algebra homomorphism $\pi$ is defined by $\pi (X\otimes P)=
X\otimes P$ and $\pi (K)=0.$
 \vskip 3mm
 \flushpar
 {\bf (A.2)  Irreducible representations of}
 $\tilde{\frak{g}}$.  Fix an irreducible  (finite dimensional) representation
$V$
   of $\frak{g}$ and a
 number $\ell \in \Bbb C$ (to be called the {\it central charge}).  Then we
define the associated {\it generalized Verma
 module} for $\tilde{\frak{g}}$ as
 $$M(V,  \ell)= U (\tilde{\frak{g}})
 \otimes_{U(\tilde{\frak{p}})} I_{\ell} (V),$$
 where the $\tilde{\frak{p}}$-module $I_{\ell} (V)$ has the underlying
 vector space same as  $V$ on which $\tilde{\frak{u}}$
 acts trivially, the central element $K$ acts via the scalar $\ell$
 and the action of $\frak{g}=\frak{g} \otimes t^o$ is
 via the $\frak{g}$-module structure on $V$.

 In the case when $\ell \neq -h$ (where $h$ is the dual Coxeter
 number of ${\frak{g}}$; cf. Remark 5.3), $M(V, \ell)$ has a unique irreducible
 quotient denoted $L(V, \ell)$. We assume in the sequel that $\ell \neq -h$.
(It is easy to see that the $\tilde\frak{g}$-module structure on $M(V, \ell)$
and $L(V, \ell)$ extends to a $\tilde\frak{g}_{\text{comp}}$-module structure.)

 \vskip 5mm
 \flushpar
 {\bf (A.3)  Remark.} It is easy to see that any vector $v \in
 M(V, \ell)$ is contained in a finite dimensional
 $\frak{g}$-submodule of $M(V, \ell)$.  In particular, the
 same property holds for any vector in $L(V, \ell)$.
 \vskip 3mm
 \flushpar
 {\bf (A.4)  Definition.}  Consider the Lie subalgebra
 $\frak{r}^o$ of $\tilde{\frak{g}}$ spanned by
 $\{Y_{\theta} \otimes t, \theta^{\vee} \otimes 1, X_{\theta}
 \otimes t^{-1} \}$, where $Y_{\theta} $ (resp. $X_{\theta})$ is a
 non-zero root vector of $\frak{g}$ corresponding to the root
 $-\theta$ (resp. $\theta)$ and the coroot $\theta^{\vee}$ is to
 be thought of as an element of $\frak{h}$. Then the Lie algebra $\frak r^o$ is
isomorphic with $sl(2)$ (cf. proof of Theorem 5.4).

 A $\tilde{\frak{g}}$-module $W$ is said to be {\it
 integrable} if every vector $v \in W$ is contained in a finite
 dimensional ${\frak g}$-submodule of $W$ and
 also $v$ is contained in a finite dimensional
 $\frak{r}^o$-submodule of $W$.

 The following lemma follows as a consequence of $sl
 (2)$-theory.
 \vskip1ex
 \flushpar
 {\bf (A.5)  Lemma.}  {\it The irreducible module $L(V, \ell)$
 (as in $\S$ A.2) is integrable if and only if $\ell$ is an integer
 and $\ell \geq < \lambda, \theta^{\vee} >$, where $\lambda$ is
 the highest weight of } $V$.
\vskip4ex

 \flushpar
 {\bf Appendix B. An introduction to Ind-varieties.}
\vskip3ex

 \flushpar
 {\bf (B.1) Definitions \cite{Sa}.} By an {\it ind-variety} we mean a set
 $X$ together with a filtration
 $$X_o \subseteq  X_1 \subseteq  X_2 \subseteq \dots ,$$
 such that \roster \item
  $\underset{n \geq 0}\to{\cup} X_n =X$, and \item
  Each $X_n$ is  a (finite dimensional) variety over $k$ such
 that the inclusion $X_n \hookrightarrow X_{n+1}$ is a closed
 immersion. \endroster

 An ind-variety $X$ is said to be {\it projective} (resp. {\it
 affine}) if each $X_n$ is projective (resp. affine). For an affine ind-variety
$X$, we define its {\it ring of regular functions}  $k[X]$ by $k[X] = \underset
{n \to \infty} \to{\text{Inv~lt.}}~ k[X_n].$ Putting the discrete topology on
each $k[X_n]$ and taking the inverse limit topology  on $k[X]$, we obtain
$k[X]$
 as a topological ring.

 Let $X$ and $Y$ be two ind-varieties with filtrations $X_n$ and
 $Y_n$ respectively.  A map $f:X \rightarrow Y $ is said to be a
 {\it morphism} if for every $n \geq 0$, there exists a number
 $m(n) \geq 0$ such that $f(X_n) \subseteq Y_{m(n)}$ and moreover
 $f_{\mid X_n}: X_n \rightarrow Y_{m(n)}$ is a morphism. Clearly,  a morphism
$f: X \to Y$ (between two affine ind-varieties) induces a continuous
$k$-algebra homomorphism $f^\ast : k[Y] \to k[X].$

 A morphism $f:X \rightarrow Y$ is said to be an {\it
 isomorphism} if $f$ is bijective and $f^{-1}:Y \rightarrow X$
 also is a  morphism. Two ind-variety structures on the same set $X$ are said
to be {\it equivalent} if the identity map I: $X \to X$ is an isomorphism of
ind-varieties.

 We define the {\it Zariski topology} on an ind-variety $X$ by
 declaring a set $ U \subseteq X$ open if and only if $U \cap X_n$ is
 Zariski open in $X_n$ for each $n$.
 \vskip 5mm
 \flushpar
 {\bf (B.2)  Exercises.} ~ (a)  For an ind-variety $X$ (under the
 Zariski topology), a subset $Z \subseteq X$ is closed if and only if $Z
 \cap X_n$ is closed in $X_n$ for each $n$.

 (b) A morphism $f:X \rightarrow Y$ between two ind-varieties is
 continuous.

 (c) Any continuous map $f:X \rightarrow Y$ between two ind-varieties
satisfies that for each
 $n$, there exists a  $m(n)$ such that $f(X_n) \subseteq Y_{m (n)}$.
 \vskip 5mm
 \flushpar
 {\bf (B.3) Examples.}  (1) Any (finite dimensional) variety
 $X$ is of course canonically an ind-variety, where we take each
 $X_n=X$.

 (2) If $X$ and $Y$ are ind-varieties then $X \times Y$ is
canonically an ind-variety, where we define the filtration by
 $$(X \times Y)_n =X_n \times Y_n.$$

 (3) $ {\Bbb A}^{\infty}:= \{ (a_1, a_2, a_3, \cdots)$; where all but
 finitely many $a_i'$s are  zero and each $a_i \in k\}$ is an
 ind-variety under the filtration : $\Bbb A^1 \subset \Bbb A^2 \subset \Bbb A^3
\subset
 \cdots,$ where $\Bbb A^n \subset \Bbb A^{\infty}$ is the set of
 all the sequences with $a_{n+1} = a_{n+2}= \cdots =0$, which of
 course is the $n$-dimensional affine space.

 (4) Any  vector space $V$ of countable dimension over $k$ is
 canonically an affine ind-variety: Take a basis $\{e_i\}_{i \geq
 1}$ of $V$.  This gives rise to a $k$-linear isomorphism
 $\Bbb A^{\infty} \tilde{\rightarrow} V$ (taking $(a_1, a_2, a_3, \cdots)
\mapsto \sum a_ie_i$ ).  By transporting the
 ind-variety strucutre from  $\Bbb A^{\infty}$ via this isomorphism, we
 get an (affine) ind-variety structure on $V$.  It is easy to see
 that a different choice of basis of $V$ gives an equivalent
 ind-variety structure on $V$.
 \vskip 5mm
  \flushpar
 {\bf (B.4)  Definitions.}  (a)  Let $X$ be an ind-variety with
 the filtration $(X_n)$. For any $x\in X$,  define the {\it
 Zariski tangent space} $T_{x} (X)$ of $X$ at $x$ by
 $$
 T_x(X)=\, \underset {n \to \infty} \to{\text{limit}}
 ~~T_x(X_n),
 $$
 where $T_x (X_n)$ is the Zariski tangent space of $X_n $ at $x$.
 (Observe that $x\in X_n$ for all large enough $n.)$

 A morphism $f:X \rightarrow Y$ clearly induces a linear map
 $(df)_x:T_x (X) \rightarrow T_{f(x)} (Y)$ (for any $x \in X)$,
 called the  {\it derivative} of $f$ at $x$.

 (b) An ind-variety $H$ is said to an {\it ind-algebraic group}
 (for short an ind-group), if the underlying set $H$ is a group
 such that the map
 $H \times H \rightarrow H~, $ taking $(x,y) \mapsto xy^{-1}~,$
 is a morphism.  In this article, we only have occasion to consider affine
 ind-groups, i.e., ind-algebraic groups $H$ such that  $H$ is an affine
ind-variety.

By a {\it group morphism} between two ind-groups $H$ and $K$, we mean a group
homomorphism $f : H \to K $ such that $f$ is also a morphism of ind-varieties.

 An abstract representation of the ind-group $H$ in a countable
 dimensional $k$-vector space $V$ is said to be {\it algebraic} if
 the map $H \times V \rightarrow V , $ defined by  $(h,v) \mapsto h.v$,  is a
 morphism.

For an ind-group $H$ and ind-variety $Y$, we say that $Y$ is an $H$-{\it
variety} if the group $H$ acts on $Y$ such that the action $H\times Y \to Y$ is
a morphism of ind-varieties.
 \vskip 3mm
 \flushpar
 {\bf (B.5)  Proposition \cite{Sa}.} {\it For an ind-group $H$, the
 Zariski tangent space $T_e (H)$ at the identity element $e$ is  endowed with a
 natural Lie algebra structure.  We denote this Lie algebra by
Lie $H$.

Moreover, if $\alpha : H \to K $ is a group morphism between two ind-groups,
then the induced map $(d\alpha)_e : \text{Lie}~ H \to \text{Lie}~ K $ is a Lie
algebra homomorphism.}
\vskip 3mm
\flushpar
{\bf Proof}.
 Denote $k[H]$ by $A$.  The multiplication map $\mu
 =\mu _H: H\times H\to H$ taking $(h_1,h_2)\mapsto
 h_1h_2$ induces a continuous homomorphism $\mu^\ast :
 A\to k[H\times H]$.  There is a canonical inclusion
 $A\otimes A\hookrightarrow k[H\times H]$, and it is easy
 to see that the image is dense in $k[H\times H]$.  So we
 denote $k[H\times H]$ by $A\hat{\otimes}A$, and view it
 as a certain completion of $A\otimes A$.  Let $\epsilon
 :A\to k$ be the homomorphism, taking $f\mapsto f(e)$.
 Let $\frak{m}= \ker \epsilon $.  Then for any $f\in
 \frak{m}$
 $$
 \mu^\ast f - f\otimes 1 - 1\otimes f\in
 \frak{m}\hat{\otimes}
 \frak{m},
 \tag1
 $$
 where $\frak{m}\hat{\otimes} \frak{m}$ denotes the
 closure of $\frak{m} \otimes \frak{m}$ in $A\hat{\otimes
 }A$.

 A continuous derivation $D: A\to A$ is said to be {\it
 invariant} if $L_h^\ast \circ D = D\circ L^\ast _h$, for
 all $h\in H$, where $L^\ast _h:A\to A$ is the algebra
 homomorphism induced from the left translation map
 $L_h:H\to H$ taking $g\mapsto hg$.  The set Der$\,A$
 of continuous invariant derivations of $A$ is a Lie
 algebra  under
 $$
 [D_1,D_2] := D_1\circ D_2 - D_2\circ D_1\, ,\qquad
 D_1,D_2 \in \text{Der}\,A.
 $$

 Define the map $\eta: T_e(H)\to \,\text{Der}\,A$ as
 follows.  Take $v\in T_e(H)$.  Then $v\in T_e(H_n)$ for
 some $n$ (where $H_n$ is the filtration of $H$).  By
 definition, $T_e(H_n) =
 \,\text{Hom}\,_k(\frak{m}_n/\frak{m}^2_n,k)$, where
 $\frak{m}_n:= \{ f\in k[H_n] : f(e)= 0\}$ is the maximal
 ideal of $k[H_n]$ corresponding to the point $e$.  In
 particular, $v$ gives rise to a $k$-linear map $\hat{v}:
 \frak {m}_n\to k$.  Let $\bar{v}: A\to k$ be the
 continuous linear map such that $\bar{v}(1)=0$, and
 $\bar{v}_{\vert_{\frak {m}}} = \hat{v}\circ \pi_n$,
 where $\pi_n:\frak {m}\to \frak {m}_n$ is the canonical
 restriction map.  Now the map $\eta (v): A\to A$  is
 defined by
 $$
 \eta (v) =
 ( I \hat{\otimes} \bar{v}) \circ \mu^\ast,
 $$
 where $I:A\to A$ is the identity map and $I
 \hat{\otimes}\bar{v} : A\hat{\otimes}A \to A\hat{\otimes
 }k=A$ is the completion of the map $I\otimes
 \bar{v}$.  By using (1), we get that $\eta(v)$ is a
 derivation.  Further, it can be seen that $\eta(v)$ is
 invariant and hence $\eta(v)\in $ Der$\, A$.

 Conversely, we define a map $\xi: \,\text{Der}\, A\to
 T_e(H)$ as follows.  Take $D\in $ Der$\,A$ and consider
 $\epsilon\circ D_{\vert_{\frak {m}}}: \frak{m}\to k$.  Since
 $D$ is continuous, there exists some $n$ such that
 $\epsilon\circ D_{\vert_{\frak {m}}}$ factors through
 $\frak{m}_n$, giving rise to a map (denoted ) $\beta
 _D:\frak{m}_n\to k$.  Since $D$ is a derivation, $\beta
 _D(\frak{m}_n^2)=0$ and hence $\beta _D$ gives rise to
 an element $\hat{\beta }_D \in T_e(H_n)$.  Now set
 $\xi(D)= \hat{\beta }_D$.

 It can be easily seen that $\xi \circ \eta$ and
 $\eta\circ \xi$  both are the identity maps, in
 particular, $\eta$ and $\xi $ are isomorphisms.  We
 now transport the Lie algebra structure from Der$\,A$ to
 $T_e(H)$ (via $\eta$).

 Finally, we prove that for any group
morphism $\alpha
 :H\to K$, the induced map $\dot{\alpha}= (d\alpha
 )_e:\,\text{Lie}\, H\to \,\text{Lie}\, K$ is a Lie
 algebra homomorphism:

 To prove this, it suffices to show that the following
 diagram is commutative (for any $v\in \,\text{Lie}\,
 H$):
$$\CD
 k[K]  @>\alpha^\ast>>  k[H]\\
 @V\eta(\dot{\alpha}v)VV   @V\eta(v)VV\\
 k[K]  @>\alpha^\ast>>  k[H].
 \endCD\tag"($\ast$)"$$

 Take $f\in \frak{m}_K$, where $\frak{m}_K\subset k[K]$
 is the maximal ideal corresponding to the point $e$.
 Then, by the definition of the map $\eta$,
 $$\align
 \eta(v) (\alpha^\ast f)
 &=(I \hat{\otimes} \bar{v})
 \mu^\ast _H (\alpha^\ast f), \qquad \text{and}
 \tag2\\
 \eta (\dot{\alpha}v) f
 &=(I \hat{\otimes} (\overline{\dot{\alpha}v}))
 \mu^\ast _K (f).\tag3\\
 \intertext{Further,}
 \alpha^\ast \eta(\dot{\alpha}v)f &=
 (\alpha^\ast \hat{\otimes} (\overline{\dot{\alpha}v}))
 \mu^\ast _K(f),\qquad \text{whereas} \tag4\\
  \overline{\dot{\alpha}v} &= \bar{v}\circ \alpha
^\ast, \qquad \text{and} \tag5\\
 (\alpha^\ast \hat{\otimes}\alpha^\ast) \circ\mu
^\ast _K &=
 \mu^\ast _H\circ\alpha^\ast .\tag6
 \endalign$$
 Now combining (2)--(6), we get the commutativity of the
 diagram ($\ast$).  This proves the proposition.\qed
 \vskip1ex

 \flushpar
 {\bf (B.6)   Lemma.}  {\it An algebraic representation $\theta$ of
 an ind-group $H$ in a (countable dimensional) vector space $V$
 induces (on `differentiation' as defined below) a representation  $d \theta$
of
 the Lie algebra Lie $(H)$ on the same space $V$.}
 \vskip 3mm
 \flushpar
 {\bf Proof}.  Fix $ v \in V$ and consider the map
 $\theta_v:H \rightarrow V ,$ taking $ h \mapsto h v.$
 Consider the derivative
 $(d \theta_v)_e: T_e (H) =$ Lie $(H) \rightarrow T_v (V) \approx
 V.$
 Then the representation $d \theta$ : Lie $(H) \times V
 \rightarrow V$ is defined as $(x,v) \mapsto (d \theta_v)_e (x)$.
  We claim that $d \theta$ is a Lie algebra representation:

 We abbreviate $d\theta(x,v)$ by $x.v ~.$ For any $v\in V$,
 define the evaluation map $e(v): k[V]\to k $ by $ e(v)f =
 f(v)$, for $f\in k[V]$.  Fix any $v_0\in V$.  Then
 $v\in T_{v_0}(V) \approx V$  can be thought of (by the definition of the
 tangent vector) as a $k$-linear map
 $\bar{v}: k [V]\to k$, such that $\bar{v}(1)=0$. If
 $v,w\in T_{v_0}(V)$ are such that
 $\bar{v}_{\vert_{V^\ast}} = \bar{w}_{\vert_{V^\ast}}$,
 then $\bar{v}=\bar{w}$, where $V^\ast \subset k[V]$
 denotes the full vector space dual of $V$.  Moreover, as
 is easy to see,
 $$\gather
 \bar{v}_{\vert_{V^\ast}} = e
(v)_{\vert_{V^\ast}} \, .
 \tag1\\
 \intertext{By definition (for any $v\in V$ and $x\in
 T_e(H)$)}
 \overline{x\cdot v} = (\bar{x}\hat{\otimes}e(v))\circ
 \theta^\ast , \tag2\\
 \intertext{where $\theta:H\times V\to V$ is the
 representation.  Since $\theta$ is linear in the
 $V$-variable,}
 \theta^\ast (V^\ast)\subset k[H] \hat{\otimes}V^\ast
 .\tag3\\
 \endgather$$

 Consider the following commutative diagram, where $A=
 k[H]$, and $I$ stands for the identity maps.
 $$
   \gather
\gathered k[V]\\  \sideset\and{\theta^\ast}\to\downarrow
 \\A\hat{\otimes} k[V] \endgathered \\
 \vspace{1\jot}
  \sideset_{\mu^\ast\hat{\otimes}I}\and\to\swarrow
  \qquad\qquad\qquad
  \sideset\and_{I\otimes\theta^\ast}\to\searrow\\
 \vspace{2\jot}
  (A\hat{\otimes}A) \hat{\otimes} k[V]  \quad \approx
 \quad A\hat{\otimes}(A\hat{\otimes}k[V])\\ \vspace{2\jot}
  \sideset_{(\bar{x}\hat{\otimes}\bar{y})\hat{\otimes
 }e(v)}\and\to\searrow
  \qquad\qquad\qquad
  \sideset\and_{\bar{x}\hat{\otimes} (\bar{y}\hat{\otimes
 }e(v))}\and\to\swarrow\\
 \vspace{2\jot}
 k\otimes k\otimes k=k.
   \endgather
 $$

 The commutativity of the above diagram and (1)--(3)
 give the following (for any  $x,y\in T_e(H)$ and $v\in
 V$).
 $$\gather
 e(x\cdot (y\cdot v))_{\vert_{V^\ast}} =
 (((\bar{x} \hat{\otimes}\bar{y}) \hat{\otimes}e(v))
 \circ
 (\mu^\ast \hat{\otimes}I)\circ\theta^\ast
 )_{\vert_{V^\ast}}.
 \tag4\\
 \intertext{By (4) we get}
 (e(x\cdot (y\cdot v)) - e(y\cdot (x\cdot
 v)))_{\vert_{V^\ast}}
 =
 (((\bar{x} \hat{\otimes}\bar{y}-\bar{y}
 \hat{\otimes}\bar{x})
 \hat{\otimes}e(v)) \circ
 (\mu^\ast\hat{\otimes}I ) \circ \theta^\ast
 )_{\vert_{V^\ast}}.\tag5\\
 \intertext{But, as can be easily seen from the
 definition of the bracket in $T_e(H)$ (cf. proof of
 Proposition B.5),}
 (\bar{x} \hat{\otimes}\bar{y} - \bar{y}\hat{\otimes
 }\bar{x})
 \circ \mu^\ast = \overline{[x,y]} .
 \tag6
 \endgather$$
 In particular,
 $$\align
 e(x\cdot (y\cdot v) -y\cdot (x\cdot v))_{\vert_{V^\ast
 }}
 &= (( \overline{[x,y]}\hat{\otimes
 }e(v))\circ\theta^\ast)_{\vert_{V^\ast}},\\
 &= ( \overline{[x,y]\cdot v} )_{\vert_{V^\ast}}, \
 \text{by (2)} \\
 &= e ([x,y]\cdot v)_{\vert_{V^\ast}}, \ \text{by
 (1)}.
 \endalign$$

 This gives that $x\cdot (y\cdot v) - y\cdot (x\cdot v) =
 [x,y]\cdot v,$ proving the lemma.\qed

  \vskip 3mm
 \flushpar
 {\bf (B.7)  Definitions.}
For any ind-variety $Y$, by an {\it algebraic vector bundle  of rank } $r$ over
$Y$, we mean an ind-variety $E$ together with a morphism $\theta : E \to Y$
such that  (for any $n$) $E_n \to Y_n $ is an algebraic vector bundle of rank
$r$ over the (finite dimensional) variety $Y_n$ ,  where $\{ Y_n\}$  is the
filtration of $Y$ giving the ind-variety structure  and $E_n := \theta ^{-1}
(Y_n)$ . If $r=1$ , we call $E$ an {\it algebraic line bundle} over $Y$.

Let $E$ and $F$ be two algebraic vector bundles over  $Y$. Then a  morphism (of
ind-varieties) $\varphi : E \to F$ is called
a {\it bundle morphism} if the following diagram is commutative :
   $$
       \gather
       E
       \overset\varphi\to{\longrightarrow}
    F \\
       \searrow  \qquad\swarrow\\
       Y
       \endgather
       $$
and moreover $\varphi_{|E_n} : E_n \to F_n $ is a bundle morphism for all $n$.
In particular, we have the notion of isomorphism of vector bundles over  $Y$.

We define Pic
  $Y$ as the set of isomorphism classes of algebraic line bundles on $Y$. It is
clearly an abelian group under the tensor product of line bundles.

We similarly define the notion of principal $G$-bundles on an ind-variety (for
a finite dimensional algebraic group $G$).

For an ind-group $H$ and $H$-variety $Y$ (cf. \S B.4),
 an algebraic vector bundle $\theta:E\to Y$ is said to be
 an $H$-{\it equivariant vector bundle} if the
 ind-variety $E$ is also an $H$-variety such that the
 following diagram is commutative:
 $$\CD
 H \times E  @>>>  E\\
 @VV{I\times\theta}V  @VV{\theta}V\\
 H\times Y  @>>>  Y~,
 \endCD\ \ $$
  and moreover for any $y\in Y$ and $h\in H$ the fiber
 map $h \times \theta^{-1} y\to \theta^{-1}(hy)$ is linear.

\vskip4ex

 \flushpar
 {\bf Appendix C. Affine Kac-Moody groups and their flag varieties.}

\vskip3ex

Let $ G$ be a connected simply-connected simple algebraic group and let $\Cal G
:= G({\Bbb C}((t))), \Cal P :=  G({\Bbb C}[[t]])$.
 We fix a Borel subgroup $ B\subset G$ and a
 maximal torus $T\subset B$, and define the {\it standard
 Borel subgroup} $\Cal B$ of $\Cal G$ as $ev_0^{-1}(B)$,
 where $ev_0:\Cal P=G({\Bbb C}[[t]])\to G$ is the
 group homomorphism induced from the ${\Bbb C}$-algebra
 homomorphism $ {\Bbb C}[[t]]\to{\Bbb C}$, taking
 $t \mapsto 0$.

Let $N(T)$ be the normalizer of $T$ in $G$ and consider the set Mor $({\Bbb
C}^\ast, N(T))$ of all the regular maps $f : {\Bbb C}^\ast \to N(T)$, which is
a group under poinwise multiplication. Then $T$ can be thought of as a (normal)
subgroup of Mor $({\Bbb C}^\ast, N(T))$ consisting of constant loops in $T$.
Then the  {\it affine Weyl group}   $\widetilde{W}$ of $G$ is by definition
$\widetilde{W}=$ Mor$({\Bbb C}^\ast, N(T))/T$.   Clearly the (finite) Weyl
group   $W :=N(T)/T$  of $G$ is a subgroup of $\widetilde{W}.$
 \vskip2ex
 \flushpar {\bf (C.1) Bruhat Decomposition.} We can view
Mor$({\Bbb C}^\ast, N(T))$ as a subgroup of $\Cal G$. In particular, any
element $w\in $ Mor$({\Bbb C}^\ast, N(T))$ can be thought of as an element
(denoted by the same symbol) $w$ of $\Cal G$.  The
  generalised flag variety $ X:= {\Cal G}/{\Cal P} $ has the
 following {\it Bruhat decomposition}:
 $$
 X=\bigcup_{\frak w\in\widetilde{W}/W}\Cal B\frak w\Cal
 P/\Cal P\,,
 \tag"(1)"$$
where the notation $\Cal B\frak w\Cal
 P/\Cal P$ means $\Cal B w\Cal
 P/\Cal P$ for any choice of the coset representative $w$ of $\frak w$. (The
set $\Cal B\frak w\Cal
 P/\Cal P$ is independent of these choices.)
 Moreover the union in (1)
 is disjoint.

 The affine Weyl group $\widetilde{W}$ is a Coxeter group and hence has a
 Bruhat partial order $\leq $.  This induces a partial
 order (again denoted by) $\leq $ in $\widetilde{W}/W$
 defined by
 $$
 \frak u := u \text{ mod }W \leq \frak v\qquad \text{(for
 $u,v\in \widetilde{W}$)}
 $$
 if and only if there exists a $w\in W$ such that
 $$
 u\leq vw\,.
 $$
 \vskip1ex
 \par We define the {\it generalised Schubert variety}
 $X_{\frak w}$ (for any $\frak w\in \widetilde{W}/W$) by
 $$
 X_{\frak w} := \bigcup_{\frak v\leq \frak w} \Cal
 B\frak v\Cal P/\Cal P\,.
 \tag"(2)"$$

 Then clearly $X_{\frak v}\subseteq X_{\frak w}$ if and
 only if $\frak v\leq \frak w$.

\vskip1ex

 \flushpar {\bf (C.2)
Definition.} Let $\frak g$ be the
 Lie algebra of $G$. We define the {\it adjoint
 representation} Ad of $\Cal G$ in $\widetilde{\frak
 g}_{\text{comp}}$ as follows (cf. \cite{PS, Proposition
 4.3.2}):
\vskip1ex

  Embed $G\hookrightarrow SL_N$ and define for any
 $g\in\Cal G$, $Y\in\frak g\otimes\Bbb C((t))$ and
 $z\in\Bbb C$
 $$
 \text{Ad}(g)(Y+zK)=
 gYg^{-1}+
 \biggl( z-\underset{t=0}\to{\text{Res}}
 \bigg\langle g^{-1}\frac{dg}{dt},Y\bigg\rangle_{\negmedspace\negthickspace t}
 \biggr)K,
 $$
 where $\langle \,,\, \rangle_{\negthinspace t}$ is the $\Bbb
 C((t))$-bilinear extension of the normalized Killing
 form $\langle\, ,\, \rangle$ on $\frak g$ (cf. \S A.1).

\vskip2ex

 The following lemma is well known, but we give a proof
 for completeness.  Even though we do not need, a more general lemma (where the
base
 field $\Bbb C$ is replaced by any $\Bbb C$-algebra) for
 $G=SL_N$ is proved by Faltings (cf. \cite{BL, Appendix
 A}).

 \proclaim{(C.3) Lemma} Let $\pi:\tilde{\frak g}\to\text{{\rm
 End}}\,W
$ be an integrable highest weight (in particular irreducible) representation
 of $\tilde{\frak g}$.  Then there exists a unique group
 homomorphism $\hat{\pi}:\Cal G\to  PGL(W)$,
 such that the following holds for any $g\in \Cal G$ and $X\in
 \tilde{\frak g}_{\text{comp}}$ :
 $$
 \hat{\pi} (g) \bar{\pi} (X)
 \hat{\pi} (g)^{-1} =  \bar{\pi} (\text{{\rm Ad}}\,g \,X),
 \tag"($\ast$)"$$
 where $\bar{\pi }:\tilde{\frak g}_{\text{{\rm comp}}}\to
 \text{{\rm End}}\, W$ is the extension of $\,\pi$ ({\rm cf.
 \S A.2}), $PGL(W):=GL(W)/{\Bbb C^\ast}$, and $GL(W)$ is
 the group of all linear automorphisms of $W$. (We view
 $\hat{\pi}(g)\bar{\pi }(X)\hat{\pi}(g)^{-1}$ as an
 element of {\rm End}$\,$ W by taking any lift of
 $\hat{\pi}(g)$ in GL(W).) \endproclaim

 \flushpar {\it Proof.}  Fix $g\in \Cal G$. We first
 prove that if there exists an element $\theta\in
PGL(W)$
 such that $\theta\bar{\pi}(X)\theta^{-1}=\bar{\pi}
 (\text{Ad}\,g\,X)$, for all $X\in \tilde{\frak
 g}_{\text{comp}}$, then $\theta$ is unique:

 For, if possible, let $\theta_1$ be another such element.  Then
 $$
 (\theta_1^{-1}\theta)
 \bar{\pi}(X) (\theta_1^{-1}\theta)^{-1}
 = \bar{\pi}(X) ,\qquad \text{for all }X\in\tilde{\frak
 g}_{\text{comp}}.
$$
 But $W$ being irreducible, $(\theta_1^{-1}\theta)=1$ (in
 $PGL(W)$).  This proves the uniqueness assertion.
 \vskip1ex

 Define the set
 $$
 S=\{ g\in\Cal G:\hat{\pi}(g)\text{ is defined
 satisfying ($\ast$) for all }X\in \tilde{\frak
 g}_{\text{comp}}\}~.
 $$
 By uniqueness, it is clear that $S$ is a subgroup of $\Cal G$
 and moreover the map $\hat{\pi}: S\to PGL(W)$ is a group
 homomorphism. We next prove that $S=\Cal G$ :

 For any root vector $X\in \frak g$ and $p\in \Bbb
 C((t))$, define
 $$
 \hat{\pi}( \exp(X\otimes p)) = \exp (\bar{\pi}(X\otimes
 p)).
 \tag1$$
 (Since $X$ is a root vector and $W$ is integrable,
 $\bar{\pi}(X\otimes p)\in \text{End}\,W$ is locally
 nilpotent, in particular, $\exp(\bar{\pi}(X\otimes p))$
 is well defined.)  It is easy to see that
  $\hat{\pi}(\exp(X\otimes p))$, as defined by (1),  satisfies
 $(\ast)$ for every $X\in \tilde{\frak g}_{\text{comp}}$.
 Further, by a result of Steinberg, the group
 generated by the elements $\exp(\bar{\pi}(X\otimes p))$
 is the whole group $\Cal G$.  This proves the
 lemma.\qed
\vskip2ex

 \subhead (C.4) Central Extension \endsubhead  Recall the
 definition of the integrable highest weight
 (irreducible) $\tilde{\frak g}$-module $W_0=L(\Bbb C,1)$
 with central charge 1 from Lemma (A.5), where $\Bbb C$
 is the trivial one-dimensional representation of $\frak
 g$.  By Lemma (C.3), there exists a group homomorphism
 $\hat{\pi}: \Cal G \to PGL(W_0)$.  We define the group
 $\tilde{\Cal G}$ as the pull-back $\hat{\pi}^\ast(GL(W_0))$:
 $$\CD
 \tilde{\Cal G}  @>>>  GL(W_0)\\
 @VVV  @VVV\\
 \Cal G @>>\hat{\pi}>  PGL(W_0)\ .
 \endCD
 $$
 Then $\tilde{\Cal G}$ is a central extension:
 $$
 1 \to \Bbb C^\ast \to \tilde{\Cal G}\overset{\pi}\to\longrightarrow\Cal G\to
1.\tag1
$$
 By the very definition of $\tilde{\Cal G}$, the $\tilde{\frak
 g}$-module $W_0$ becomes a $\tilde{\Cal G}$-module.  In
 particular, the tensor product $W_0^{\otimes m}$ (for
 $m\geq 0$) acquires a canonical $\tilde{\Cal G}$-module
 structure.  Since the integrable $\tilde{\frak
 g}$-module $L(\Bbb C,m)$ with central charge $m\geq 0$
 is a $\tilde{\frak g}$-submodule of $W_0^{\otimes
 m}$, it
 is easy to see (using (1) of the proof of Lemma C.3)
 that $\tilde{\Cal G}$ keeps $L(\Bbb C,m)\subset W_0^{\otimes
 m}$ stable.  In particular, $L(\Bbb C,m)$ acquires a
 canonical $\tilde{\Cal G}$-module structure.
\vskip1ex

  \flushpar
  {\bf (C.5)} {\bf Realizing $X$  as an ind-variety via
 representation theory}.
 Define the filtration $\{X_n\}_{n \geq 0}$ of $X$ as follows:

 $$X_n:= \cup  X_{\frak{w}},$$
 where the union is taken over those $\frak{w}=w$ mod $W \in
 \widetilde{W} / W$ such that $\ell (w) \leq n$.

 Fix any integer $\ell > 0$ and consider the irreducible
 $\tilde{\frak{g}}-$module $L(\Bbb C, \ell) $.  From the Bruhat decomposition,
it is easy
 to see that the map $i=i_\ell: X \rightarrow \Bbb P (L (\Bbb C, \ell)), g
 {\Cal P} \mapsto [gv_o]$  (where $v_o$ is the highest weight
 vector of $L(\Bbb C, \ell)$ and $\Bbb P (L (\Bbb C, \ell))$ denotes the space
of lines in $L (\Bbb C, \ell)$) is injective.
 As proved in \cite{Sl, \S2.4}, for any $n$, there exists a finite
 dimensional subspace $W_n \subset L (\Bbb C, \ell)$ such that
 $i(X_n) \subset \Bbb P (W_n)$ and moreover $i(X_n)$ is Zariski
 closed in $\Bbb P (W_n)$.  We endow $X_n$ with the projective (reduced)
variety structure so that $i_{\mid X_n}$ is a closed immersion.
 This makes $X$ into a projective ind-variety.  Further, the
 ind-variety structure does not depend upon  the particular choice
 of $\ell >0 $ (cf. \cite{Sl, \S2.5}). Equipped with this ind-variety
structure, we denote $X$ by $X_{\text{rep}}$ .

For any ${\frak w}  \in
 \widetilde{W} / W$ the generalised Schubert variety $X_{\frak w}$ is an
irreducible Zariski closed subset of $X_{\text{rep}}.$ We endow $X_{\frak w}$
with the projective (reduced) variety structure so that $X_{\frak w}
\hookrightarrow X$ is a closed immersion. Then $\Cal
 B\frak w\Cal P/\Cal P \subset X_{\frak w}$ is an open subset which is
biregular isomorphic with the affine space $\Bbb
 C^{\ell(\frak w)}$, where   $\ell(\frak w)$ is the length of
 the smallest element in the coset $\frak w$.

Since $X_n$ is a variety, we can equip $X_n$ with the (Hausdorff) analytic
topology and put the inductive limit topology on $X$. The decomposition (1) of
\S (C.1) provides a cellular decomosition of $X$, making it into a CW complex.
 \vskip2ex

 \vskip1ex
Following \cite{Sl,\S
 2.7}, we define the homogeneous line bundles on
 $X_{\text{rep}}$:
\vskip1ex
 \flushpar {\bf (C.6) Definition.}
 For any countable dimensional vector space $V$, we first
 define the tautological line bundle $\Cal L_V$ on $\Bbb
 P(V)$ as follows:  Consider the subset
 $$
 \Cal L_V = \{ (\ell,v) \in \Bbb P (V)\times V: v\in
 \ell\} .
 $$
 Then $\Cal L_V$ is a Zariski-closed subset of the
 ind-variety $\Bbb P(V)\times V$.  We equip $\Cal L_V$
 with the ind-variety structure so that $\Cal
 L_V\hookrightarrow \Bbb P(V)\times V$ is a closed
 immersion.  Now the projection on the first factor $\Cal
 L_V\rightarrow \Bbb P(V)$ realizes $\Cal L_V$ as an
 algebraic line bundle on $\Bbb P(V)$.

 For any $\ell >0$, define the algebraic line bundle
 $\frak L(\ell\chi_0)$ on $X$ as the pull-back of the
 dual $\Cal L^\ast$ of the tautological line bundle $\Cal
 L= \Cal L_{L(\Bbb C,\ell)}$ on $\Bbb P(L(\Bbb C,\ell))$
 via the embedding $i_\ell :X\to \Bbb P(L(\Bbb C,\ell))$
 of the above section.  For any integer $\ell <0$, we
 define the line bundle $\frak L(\ell\chi_0)$ as the
 dual $\frak L(-\ell\chi_0)^\ast$ and for $\ell =0$,
 $\frak L(\ell\chi_0)$ is defined to be the trivial
line
 bundle.  It is easy to see (cf. \cite{Sl, \S 2.7}) that
 the line bundle $\frak L(\ell\chi_0)$ is isomorphic
 with the line bundle $\frak L(\chi_0)^{\otimes \ell}$.

 The group ${\Cal G}$ acts (set theoretically) on
 $X=\Cal G/\Cal P$ via $g(h\Cal P) = gh\Cal P$, for
 $g,h \in {\Cal G}$.   We
 denote the action of $g\in {\Cal G}$ on $X$ by
 $L_g$.  This action lifts to an action of the group
 $\tilde{\Cal G}$ on the line bundle $\frak L(\ell\chi
 _0)$ (for $\ell <0$) via
 $$
 g(x,v) = (L_{\pi(g)}x,gv),\qquad \text{for any }g\in \tilde{\Cal
 G},\ x\in X \text{ and }v\in i_\ell (x),
 $$
 where
$\pi:
 \tilde{\Cal G} \to \Cal G$ is the canonical map.
 Observe that for any fixed $g\in \tilde{\Cal G}$, the
 action of $g$ on $\frak L(\ell\chi_0)$ is an algebraic
 automorphism of the algebraic line bundle $\Cal
 L(\ell\chi_0)$ inducing the  automorphism
 $L_{\pi(g)}$ on the base $X$.

 Set $\frak{L}_o := \frak{L}(-\chi_0)$, and define the   Mumford group (cf.
\cite{PS, Remark (i), page 115})
 $
 \text{Aut}(\frak{L}_o ) =
 \{ (g,\varphi) :g\in \Cal G$ and $\varphi$
 is an algebraic automorphism of the line bundle
 $\frak{L}_o$
 inducing the map $L_g$ on the base$\}$.
  Then Aut$(\frak{L}_o )$ is a group under
 $$
 (g_1,\varphi_1) (g_2,\varphi_2) = (
 g_1g_2,\varphi_1\varphi_2) .
 $$
 The projection on the first factor gives a group
 homomorphism $\delta :\, \text{Aut}\, (\frak{L}_o
 )\to \Cal G$.  Since $\tilde{\Cal G}$ acts on
 $\frak{L}_o$, there is a canonical group
 homomorphism $\xi : \tilde{\Cal G}\to \,\text{Aut}\,
 ({\frak L}_o)$ making the following diagram
 commutative:
  $$
   \matrix\format\c&&\quad\c\\
   \tilde{\Cal G}  &\overset{\xi}\to\longrightarrow
  &\text{Aut}\, (\frak{L}_o) \\
   \sideset{\pi}\and\to\searrow
  &&\negthickspace \sideset\and{\delta}\to\swarrow \\
   &\Cal G
   \endmatrix
  $$
 Since $\pi$ is surjective, so is $\delta $.  In
 particular, $\xi$ is an isomorphism.

\vskip1ex

We also need another  `lattice' description of the ind-variety structure on $X$
(cf. \cite{KL, \S5}).
 \vskip1ex
 \flushpar
{\bf (C.7)} {\bf Realizing $X$ as an ind-variety via lattices.} We first
consider the case of
  $G= SL_N$. Denote  $V=\Bbb C^N,$ and $ A=\Bbb C [[t]]$.
    For any $n \geq 0$, consider the set
 ${\Cal F}_n$ of  A-submodules $L \subset V \otimes_{\Bbb C} \Bbb C
 ((t))$  such that (denoting  $V \otimes_{\Bbb C} A$ by $L_o$)
$$t^n L_o
 \subset L \subset t^{-n} L_o ~, ~\text{and ~ dim~} L/t^nL_o =nN~.$$

 Let $V_n:= t^{-n} L_o/t^n L_o$ be the complex vector space of dimension $2nN$.
 Then the
 multiplication by $t$ induces a nilpotent endomorphism
 $\bar{t}_n$ of $V_n$ and hence $1+\bar{t}_n$ is a
 (unipotent) automorphism of $V_n$.  In particular,
 $1+\bar{t}_n$ induces a biregular isomorphism of the
 Grassmannian $\text{Gr}(nN,2nN)$ of $nN$-dimensional
 subspaces of the $2nN$-dimensional space $V_n$.  Let
 $\text{Gr}(nN,2nN)^{1+\bar{t}_n} $ denote its fixed point.
 Then clearly the map $j_n:\Cal F_n\to
 \text{Gr}(nN,2nN)^{1+\bar{t}_n}$ given by $L\mapsto
 L/t^n L_o$ is a bijection. We pull the
 (reduced) subvariety structure of Gr$( nN,2nN
 )^{1+\bar{t}_n}$ via $j_n$ to equip $\Cal F_n$ with a
 projective variety structure.  We next claim that the
 canonical inclusion $\Cal F_n\hookrightarrow \Cal
 F_{n+1}$ is a closed immersion:

 Consider the commutative diagram:
 $$\CD
 \Cal F_n  @>j_n>>  Gr(nN,2nN)^{1+\bar{t}_n} \\
 @VVV  @VV\theta_nV\\
 \Cal F_{n+1}  @>j_{n+1}>> ~~~~
 Gr((n+1)N,2(n+1)N)^{1+\bar{t}_{n+1}}
 \endCD
 $$
 where the map $\theta_n$ takes
 $
 W\subset t^{-n} L_o/t^n L_o \approx t^{n-1}\, V
 \oplus t^{n-2}\, V
 \oplus \dots \oplus t^{-n}\, V
 \mapsto
 t^n V\oplus W. $  It is easy to
  see that $\theta_n$
 is a closed immersion.  This equips $\Cal F= \cup_{n\geq
 0}\Cal F_n$ with a projective ind-variety structure.

 Let $\Cal G^o := SL_{N}({\Bbb C}((t)))$ and $\Cal P^o :=
 SL_N(A)$ and set $X^o= \Cal G^o/\Cal P^o$.  By
 virtue of the following lemma, the map $\beta :X^o\to
 \Cal F$ is a bijection.  By transporting the projective
 ind-variety structure from $\Cal F$ to $X^o$ (via $\beta
 $) we equip $X^o$ with a projective ind-variety
 structure.  With this structure we denote $X^o$ by
 $X^o_{\text{lat}}$.  We also define the filtration
 $\hat{X}_n^o$ of $X^o$ by
 $$
 \hat{X}_n^o = \beta^{-1}(\Cal F_n).
 $$
 The group $\Cal G^o$ acts
 canonically on $V\otimes \Bbb C((t))$.
 \proclaim{(C.8) Lemma}  The map $g\Cal P^o\mapsto
 gL_o$ (for $g\in \Cal G^o$) induces a bijection
 $\beta:X^o\to\Cal F$.
 \endproclaim

 \flushpar {\it Proof}.  Fix $g\in \Cal G^o$.  It is easy to
 see that there exists some $n$ (depending upon $g$) such
 that
 $$
 t^nL_o\subset gL_o \subset t^{-n} L_o.
 \tag1$$
 Of course $gL_o$ is $t$-stable.  We next calculate the  dimension of
  $gL_o/t^nL_o$:

 By the Bruhat decomposition (1) of \S (C.1), it suffices to
 assume that $g$ is a morphism $\Bbb C^\ast\to D$ taking
 $1\mapsto 1$, where $D$ is the diagonal subgroup of
  $SL_N$.  Write
 $$
 g(t) =
 \pmatrix t^{n_1} &&0\\ &\ddots \\0 &&t^{n_N}\endpmatrix,
 \qquad \text{for } t\in \Bbb C^\ast ~~\text{and}~n_i \in \Bbb Z.
 $$
 Then $\sum n_i=0$.  Now
 $$
 \dim (gL_o/t^nL_o) = (n-n_1) +\dots+(n-n_N)
 = Nn-\sum n_i = Nn~.
 $$
 This proves that $gL_o\in \Cal F_n$.

 Conversely, take $L\in \Cal F_n$.  Since $A$ is a PID
 and $t^kL_o$ is $A$-free of rank $N$ (for any $k\in \Bbb
 Z$), we get that $L$ is $A$-free of rank $N$.  Further,
 $L\otimes _A\Bbb C((t))\to V
 ((t))$ is an isomorphism, where $ V((t)) := V\otimes_{_{\Bbb C}}\Bbb
 C((t)).$
 Let $\{ e_1,\dots,e_N\}$ be the standard $\Bbb
 C$-basis  of $V$ and take an $A$-basis $\{ v_1,\dots
 ,v_N \}$ of $L$.  Now define the $\Bbb C((t))$-linear
 automorphism $g$ of $V((t))$ by $ge_i=v_i$
 ($1\leq i\leq N$).  We prove that $\det g$ is a
unit of
$A$:  Write $\det g=t^ku$, where $k\in \Bbb Z$ and $u$
 is a unit of $A$.  Consider the $\Bbb C((t))$-linear automorphism
$\alpha$ of $V((t))$ defined by
 $$\align
 \alpha\, e_i &=e_i,\qquad \text{for }1\leq i < N,\\ \vspace{2\jot}
 &= t^{-k} u^{-1}e_N,\qquad \text{for }i=N.
 \endalign$$
  Then det$(g\alpha) = 1,$  and
 $t^{n+\vert k \vert}L_o\subset (g \alpha)L_o\subset
 t^{-n-\vert k \vert}L_o.
 $

\noindent
 Hence, by the first part of the proof, we get
 $$
 \dim\ \biggl( \frac{g\alpha(L_o)}{t^{n+\vert k\vert}L_o}\biggr) =
 (n+\vert k \vert)N.
 \tag "(2)"$$
 On the other hand,
 $$\align
 \dim\ \biggl(\frac{g\alpha(L_o)}{t^{n+\vert k\vert}L_o}\biggr) &=
 \dim \frac{gL_o}{t^nL_o}+\vert k \vert N+k\tag "(3)" \\  \vspace{2\jot}
 &= Nn+\vert k \vert N+k\qquad \text{(since } L\in
 \Cal F_n).
 \endalign $$
 Now combining (2) and (3), we get $k=0$, hence $(
 g\alpha )L_o= gL_o=L$.  This proves the surjectivity
 of $\beta$.  The injectivity of $\beta $ is clear.  This
 proves the lemma.$\qed $

 \vskip1ex

\subhead (C.9)\endsubhead  We now come to the case of
 general (connected, simply-connected, simple) $G$.  Fix
 an embedding $G\hookrightarrow SL_N$.  This gives rise to an embedding
 $$
 X=\Cal G/\Cal P \hookrightarrow X^o = \Cal G^o/\Cal P^o.
 $$
 The filtration $\hat{X}^o_n$ of $X^o$ (given in \S C.7)
 on restriction gives the filtration $\hat{X}_n$ of $X$,
 i.e.,
 $$
 \hat{X}_n = \hat{X}_n^o \cap X.
 $$
 In (a subsequent) Lemma (C.11), we prove that $\hat{X}_n
 $ is a Zariski closed subset of $\hat{X}_n^o$. This
  allows us to put the reduced subvariety structure on
 $\hat{X}_n $ making $X$ into a projective ind-variety.
 Equipped with this ind-variety structure, we denote $X$
 by $X_{\text{lat}}$.

 \flushpar
 {\bf (C.10) Lemma}.  {\it The two filtrations $X_n$ and $\hat{X}_n$
 of $X$ are  compatible, i.e.,  for every $n$ there exists  $k(n) $
such that
 $$X_n \subseteq  \hat{X}_{k(n)} \ \text{{\rm and}} \ \hat{X}_n  \subseteq
 X_{k(n)}.$$ }
 \vskip1ex
 \flushpar
 {\bf Proof.} Fix a maximal torus $T \subset G$ and an embedding $G
 \hookrightarrow SL_N$ such that $T$ goes inside the diagonal
 subgroup $D$ of $SL_N$.  Now there is a bijection $\widetilde{W}/W
\simeq \text{ Mor}_1
(\Bbb C^*,T)$, where
 $\text{Mor}_1$ denotes the set of morphisms $\Bbb C^* \rightarrow T$ such
 that $1 \mapsto 1$.  Since the set $\{\frak w \in \widetilde{W} / W: \ell (w)
\leq
 n\}$ is  finite, it is easy to see that $X_n \subset
 \hat{X}_{k(n)}$ (for some large enough  $k(n))$.

 Conversely (for a fixed $n$), we want to show for all but finitely many $\frak
w \in
 \widetilde{W} / W$,  $({\Cal B} \frak w {\Cal P} / {\Cal P})
 \cap\hat{X}_n=\phi$:  Represent $\frak w$ as a morphism $\Bbb C^*
 \rightarrow T\hookrightarrow D$
 $$
  z \mapsto \pmatrix z^{n_1 (\frak w)}\\ &\ddots \\
  &&z^{n_ N(\frak w)} \endpmatrix ~.
 $$
 We first claim that any $\frak w$ such that $n_i(\frak w) < -n$ (for some
 $i)$ satisfies $({\Cal
 B} \frak w {\Cal P} / {\Cal P}) \cap
 \hat{X}_n =\phi$ :  If for some $b \in {\Cal B},\,b \frak w L_o \in \Cal F_n$
, then clearly ${\frak w} L_o \in  b^{-1} \Cal F_n = \Cal F_n$, a
 contradiction to the choice of $\frak w$.  Now observe that the set $\{\frak w
\in  \widetilde{W}/W : n_i(\frak w) \geq-n$ for all $i\}$ is
 finite, since $\sum n_i (\frak w) =0.$ From this,  it follows that $\hat{X}_n
\subset X_{k(n)}$ ,
 for some large enough $k(n)$. This proves the lemma. \qed
 \vskip 3mm

 \proclaim{(C.11) Lemma} With the notation as above,
 $\hat{X}_n$ is a Zariski closed subset of $\hat{X}_n^o$
 (for any $n\geq 0$).
 \endproclaim

 \flushpar {\it Proof.}  Fix $\frak w\in\widetilde{W}/W$
 ($\widetilde{W}$ is the affine Weyl group corresponding to $G$)
 and take a coset representative $w$ of $\frak w$ of
 minimal length.  Choose any reduced decomposition
 $w=s_{i_1}\dots s_{i_p}$ (where $s_j$'s are the simple
 reflections in $\widetilde{W}$), and consider the
 Bott-Samelson-Demazure-Hansen variety $Z_w$ defined in
 \cite{Sl, \S2.3}.  Let $\Cal P_j$ be the minimal
 parabolic subgroup of $\Cal G$ corresponding to the
 simple reflection $s_j$.  Recall that, set
 theoretically, $Z_w=\Cal P_{i_1}\times \dots \times\Cal
 P_{i_p}/\Cal B^p$, where $\Cal B^p$ acts on $\Cal
 P_{i_1}\times\Cal P_{i_2}\times\dots\times\Cal P_{i_p}$ from  the right via
 $$
 (x_1,\dots ,x_p)(b_1,\dots,b_p)
 = (x_1 b_1,b_1^{-1}x_2b_2,\dots,b^{-1}_{p-1}x_pb_p),
 $$
 for $x_j\in\Cal P_{i_j}$ and $b_j\in\Cal B$.

 Define the map $\theta_w:Z_w\to X$ by
 $\theta_w((x_1,\dots,x_p)\ \text{mod}\, \Cal B^p)=
 x_1\dots x_p\Cal P$.  Since Im $\theta_w= X_{\frak w}$
 (cf. \cite{Sl, \S 2.4}), by the above lemma,
 Im$(i\circ\theta_w)\subset\hat{X}^o_m$ for some $m$, where
 $i:X\hookrightarrow X^o$ is the inclusion.  By an
argument
 similar to the proof of \cite{Sl, Theorem 2.4}, it can
 be easily seen that $i\circ\theta_w: Z_w \to
 {X}^o_{\text{lat}}$ is a morphism.  In particular,
 $Z_w$ being projective, $i(X_{\frak w})$ is closed in $
 \hat{X}^o_m$.  We now prove that $i(\hat{X}_n)$ is
 closed in $\hat{X}^o_n$:

 Observe that $ \hat{X}_n$ is left $\Cal B$-stable.  Fix
 any $\frak w \in \widetilde{W}/W$ such that  $\Cal B\frak w\Cal
 P/\Cal P\subset \hat{X}_n$.  Then we claim that
 $X_{\frak w}\subset \hat{X}_n$:  There is an open
 (dense) subset $Y_w\subset Z_w$ such that $\theta_w(Y_w)
 = \Cal B\frak w\Cal P/\Cal P$.  Hence, considering the
 morphism $i\circ\theta_w: Z_w\to X^o$, we see that
 $i\circ \theta_w(Z_w)\subset \hat{X}^o_n$ (since
 $\hat{X}^o_n$ is projective).  In particular, $X_{\frak
 w} \subset \hat{X}_n$ and thus $\hat{X}_n$ is a finite
 union (by Lemma C.10) of Schubert varieties $X_{\frak w}$.  Now since
$i(X_{\frak w})$
 is closed in $\hat{X}^o_n$, so is $i(\hat{X}_n)$.
 This proves the lemma. \qed

 \vskip1ex

 \flushpar
 {\bf (C.12) Proposition.} {\it The identity map $X_{\text{rep}} \rightarrow
 X_{\text{lat}}$ is an isomorphism of ind-varieties. }
 \vskip1ex
 \flushpar
 {\bf Proof.} Embed $G \hookrightarrow $ SL$(N)$ as in \S(C.9) and follow the
same notation as in \S\S (C.7) and (C.9). By definition, $X_{\text{lat}}
\hookrightarrow
 X_{\text{lat}}^o $ is a closed immersion.  Similarly, we claim
 that $X_{\text{rep}} \hookrightarrow X_{\text{rep}}^o $ is a closed
 immersion :

 Take the  integrable highest weight module $L=L(\Bbb C,\ell)$ for
$\widetilde{\frak g^o}$
 (for any integer $\ell >0$, where $\frak g^o=sl_N$  ), and let $W\subseteq L$
be the (integrable highest
 weight)  $\widetilde{\frak g}$-module spanned by the highest weight vector of
$L$.
 Then we have
 $$
\matrix \format\c&&\quad\c\\
 X_{\text{rep}} &\hookrightarrow& X_{\text{rep}}^o
\\
 \downarrow && \downarrow \\
 \Bbb P (W) & \hookrightarrow & \Bbb P (L)
 \endmatrix
 $$
 where both the vertical maps are by definition closed immersions,
 and moreover $\Bbb P (W) \hookrightarrow \Bbb P (L)$ is of course a
 closed immersion.
 This proves that  $X_{\text{rep}} \hookrightarrow X_{\text{rep}}^o$ is a
 closed immersion.  So,  to prove the lemma, we can take $G=SL_N:$

  Fix $\frak w\in\widetilde{W}/W$ (where $\widetilde{W}$ is
 the affine Weyl group corresponding to $G= SL_N$).  By
 the proof of Lemma (C.11) (following  the same
 notation), the map $\theta_w : Z_w\to X^o
_{\text{lat}}$ is a morphism with its image precisely
 equal to $X^o_{\frak w}$.  We denote $X^o_{\frak w}$
 endowed with the reduced subvariety structure from $X^o
_{\text{rep}}$ by $X^o_{\frak w,\text{rep}}$ (and a
 similar meaning for $X^o_{\frak w, \text{lat}}$).  Then the map
 $\overline{\theta}_w :Z_w \to X^o_{\frak w,\text{rep}}$
 (the map $\overline{\theta}_w$ at the level of sets is
 nothing but $\theta_w$) is a surjective morphism (cf.
 [Sl, Theorem 2.4]) and moreover $X^o_{\frak
 w,\text{rep}}$ is a normal irreducible variety (cf.
 \cite{Ku$_1$, Theorem 2.16}).  We claim that the inclusion
 map $I_{\frak w}: X^o_{{\frak w},\text{rep}} \to
X^o_{\text{lat}}$ is a
 morphism:

 First of all, by Lemma (C.10), Image $I_{\frak w}\subset
 \hat{X}^o_n$, for some $n$.  Now the map
 $\overline{\theta}_w :Z_w \to X^o_{\frak w,\text{rep}}$
 being a proper surjective morphism, the (Zariski)
 topology on $X^o_{\frak w,\text{rep}}$ is the quotient
 topology.  Let $\Cal U \subset \hat{X}^o_n$ be an open
 subset.  Then $\theta_w^{-1}(\Cal U ) =
 (\overline{\theta}_w)^{-1} I^{-1}_{\frak w}(\Cal U )$ is
 open in $Z_w$ and hence $I^{-1}_{\frak w}(\Cal U )$ is
 open in $X^o_{\frak w,\text{rep}}$.  To prove that
 $I_{\frak w}$ is a morphism, it suffices to show that
 for any affine open $\Cal U \subset \hat{X}^o_n$, the
 map
 $
 I_{\frak w_{\vert_{I_{\frak w}^{-1}(\Cal U)}}}
 : I_{\frak w}^{-1}(\Cal U) \to \Cal U
 $
 is a morphism:   But this follows from Proposition
 (4.1), since the map $I_{\frak w} \circ
 \overline{\theta}_w=\theta_w$ is a morphism.

 Conversely, fix $n\geq 0$ and take $\hat{X}^o_n$.  Then
 $I^{-1} (\hat{X}^o_n) \subset X^o_{\text{rep}}$ is a
 closed subset and moreover (by Lemma C.10),
 $I^{-1}(\hat{X}^o_n) \subset X^o_m$ (for some $m$), in
 particular, $I^{-1} (\hat{X}^o_n) $ is a projective
 subvariety of $X^o_{\text{rep}}$.  The bijective map
 $I_n:I^{-1}(\hat{X}^o_n)\to \hat{X}^o_{n,\text{lat}}$
 (where $I_n:= I_{\vert_{I^{-1}(\hat{X}^o_n)}}$) is a
 morphism (since $X_m^o \subset X^o_\frak w$ , for some $\frak w \in
\widetilde{W}/W$).  Further, the variety
 $\hat{X}^o_{n,\text{lat}}$ is isomorphic with the
 variety $Gr(nN,2nN)^{1+\overline{t}_n}$ (cf. \S C.7).  But
 $Gr(nN,2nN)^{1+\overline{t}_n}$  is known to be
 irreducible and normal by using a result of Kostant (cf.
 \cite{Ku$_2$}).  Moreover, $I_n$ being a homeomorphism (since $I_n$ is a
proper surjective morphism),
 $I^{-1}(\hat{X}^o_n)$ is irreducible as well.  Hence by
 \cite{Mum, page 288, I. Original form}, $I_n$ is an
 isomorphism.  This shows that the identity map
 $X^o_{\text{lat}} \to X^o_{\text{rep}}$ also is a
 morphism, proving the proposition. \qed

 \vskip1ex

 So from now on, we identify $X_{\text{lat}}$ with  $X_{\text{rep}}$ and just
denote them by $X$.
 We have the following proposition determining
 Pic$\,(X)$.
 \proclaim{(C.13) Proposition} The map $\Bbb Z\to
 \text{{\rm Pic}}\,(X)$ given by
 $$
 d\mapsto \frak L(d\chi_o)
 $$
 is an isomorphism.
 \endproclaim
\demo{Proof} For  any ${\frak w}\in\widetilde{W}/W $,
 since $X_{\frak w}$ is a projective variety, by GAGA, the natural
map
$$
\text{Pic}(X_{\frak w})
\overset\sim\to\rightarrow
 \text{Pic}_{an}(X_{\frak w})
\tag"(1)"
$$
is an isomorphism, where Pic$_{an}(X_{\frak w})$ is the set of isomorphism
classes of analytic line bundles on $X_{\frak w}$.

We have the sheaf exact sequence:
$$
0 \rightarrow \Bbb Z \rightarrow \Cal O_{an} \rightarrow \Cal O^\ast_{an}
\rightarrow 0\,,
\tag"(2)"
$$
where $\Cal O_{an}$ (resp. $\Cal O^\ast_{an}$) denotes the sheaf of analytic
functions (resp. the sheaf of invertible analytic functions) on $X_{\frak w}$.
Taking the associated long exact cohomology sequence, we get
$$
\dots \rightarrow H^1 (X_{\frak w},\Cal O_{an}) \rightarrow
H^1 (X_{\frak w},\Cal O^\ast_{an}) \overset c_1\to\rightarrow
H^2 (X_{\frak w},\Bbb Z) \rightarrow
H^2 (X_{\frak w},\Cal O_{an}) \rightarrow \dots ,
\tag"(3)"
$$
where the map $c_1$ associates to any line bundle its first Chern class.  Now
$$
H^i (X_{\frak w},\Cal O) = 0,\qquad \text{for all }i > 0
\tag"(4)"
$$
by  \cite{$\text{Ku}_1$, Theorem 2.16(3)} (also proved in  \cite{M}); and by
GAGA
$$
H^i (X_{\frak w},\Cal O)
\approx H^i (X_{\frak w},\Cal O_{an}) ,
 \tag"(5)"
$$
and hence the map $c_1$ is an
isomorphism.  But
$$
 \text{Pic}_{an}(X_{\frak w}) \approx H^1 (X_{\frak w},\Cal O^\ast_{an})\,.
\tag"(6)"
$$

Hence, by combining (1) and (3)--(6), we get the isomorphism
(again denoted by)
$$
c_1:\text{Pic} (X_{\frak w}) \overset\approx\to\rightarrow
H^2(X_{\frak w},\Bbb Z) \,.
\tag"(7)"
$$

Further, the following diagram is commutative (whenever $X_{\frak w} \subseteq
X_{\frak v}$) :
$$
\CD
\text{Pic}(X_{\frak v})  @>\underset\sim\to{c_1}>>  H^2 (X_{\frak v},\Bbb Z)\\
@VVV   @VVV\\
\text{Pic}(X_{\frak w})  @>\underset\sim\to{c_1}>>  H^2 (X_{\frak w},\Bbb Z)\ ,
\endCD
\tag"$(\Cal D)$"
$$
where the vertical maps are the canonical restriction maps.  But from the
Bruhat decomposition, for any $\frak w\geq\frak s_o$, the restriction
map
$$
 H^2 (X_{\frak w},\Bbb Z) \rightarrow H^2 (X_{\frak s_o},\Bbb Z)
\tag"(8)"
$$
is an isomorphism, where $ s_o$ is the (simple) reflection corresponding to the
simple coroot $\alpha_0^{\vee}$, and $\frak s_o :=s_o $ mod $W$.  Moreover,
$X_{\frak s_o}$ being isomorphic with the
complex projective space $\Bbb P^1,\ H^2 (X_{\frak s_o},\Bbb Z)$ is a free
$\Bbb
Z$-module of rank 1, which is generated by the first Chern class $-1$ of the
line bundle $\frak L(\chi_o)_{\vert X_{\frak s_o}}$.
In particular, $\text{Pic}(X_{\frak w})$ is freely generated by $\frak
L(\chi_o)_{\vert X_{\frak w}}$, for any ${\frak w} \geq {\frak s}_o$.

 We next  prove that the canonical map $\alpha: $ Pic$~X  ~
\to
\varprojlim_ {\frak w\in\widetilde{W}/W}\,
 \text{Pic}\, (X_{\frak w})\,$ is an isomorphism:

  Since  the
 line bundles $\frak L(d\chi_0)$ (for $d \in {\Bbb Z}$) are algebraic line
bundles on $X$, the surjectivity of the map $\alpha$ follows.  Now we come to
the injectivity of $\alpha$ :

Let $\frak L \in $ Ker $\alpha$.  Fix a non-zero vector $v_o$ in the fiber of
$\frak L$ over the base point $\frak e \in X $.  Then
${\frak L}_{\vert  X_{\frak w}} $ being  a trivial line bundle on each
$X_{\frak w}$,
 we can choose a  nowhere-vanishing section
$s_{\frak w}$
of  ${\frak L}_{\vert X_{\frak w}} $ such that  $s_{\frak w} (\frak e) =  v_o
.$
We next show that for any $\frak v \geq \frak w ,  s_{\frak v_{\vert X_{\frak
w}}}=  s_{\frak w} :$  Clearly  $ s_{\frak v_{\vert X_{\frak w}}} =  f s_{\frak
w} ,$ for some regular function $f : X_{\frak w} \to
{\Bbb C}^\ast$ . But $X_{\frak w}$ being projective and irreducible,
$f$ is constant and in fact $f\equiv 1$ since $s_{\frak v} (\frak e) = s_{\frak
w} (\frak e) .$ This gives rise to a nowhere-vanishing (regular) section
$s$ of $\frak L$ on the whole of $X$ such that  $s_{\vert X_{ \frak w}}=
s_{\frak w} $. From this it is easy to see that $\frak L$ is isomorphic with
the trivial line bundle on $X$. This proves that $\alpha$ is injective, thereby
completing the proof of the proposition. \qed
\enddemo

\vskip6ex

\noindent
Deptt. of Mathematics, University of North Carolina, Chapel Hill, N. C.
27599-3250, USA.

\enddocument